\documentclass[11pt]{article}
\usepackage{pxfonts}
\usepackage{yfonts}
\usepackage{dsfont}
\usepackage{graphicx}
\usepackage{relsize}

\def\bel{\begin{equation}\label}
\def\eeq{\end{equation}}
\def\ds{\displaystyle}
\def\endproof{\hphantom{MM}
\hfill\llap{$\square$}\goodbreak}

\def\mt{\longrightarrow}
\def\v{\vskip 1em}
\def\vsk{\vskip 40em}
\def\ve{\varepsilon}
\def\R{\mathbb R}

\def\C{\mathfrak{C}}

\def\N{{\bf N}}

\def\S{\mathcal{S}}

\def\Q{{\bf Q}}

\def\A{{\bf A}}

\def\L{{\bf L}}
\def\U{\mathcal{U}}
\def\V{\mathcal{V}}

\def\setJ{\mathcal{J}}
\def\setI{\mathcal{I}}

\def\I{{\bf I}}

\def\alpha{\alphaup}
\def\beta{\betaup}
\def\gamma{\gammaup}
\def\delta{\deltaup}
\def\xi{{\xiup}}
\def\eta{{\etaup}}
\def\tau{{\tauup}}
\def\rho{{\rhoup}}
\def\phi{{\phiup}}
\def\psi{{\psiup}}
\def\lambda{{\lambdaup}}
\def\omega{\omegaup}
\def\varphi{{\varphiup}}
\def\gamma{{\gammaup}}

\def\t{{\bf t}}

\newtheorem{prop}{Proposition}[section]
\newtheorem{remark}{Remark}[section]

\begin{document}
 \[\begin{array}{cc}\hbox{\LARGE{\bf Stein-Weiss inequality on product spaces}}
 \end{array}\]
 
\[\hbox{Zipeng Wang}\]

 \begin{abstract}
We give the classification between   weighted norm inequalities of strong fractional integral operators  and their associated multi parameter Muckenhoupt characteristics,  by considering  the weights to be power functions. As a result, we extend the classical Stein-Weiss theorem  to product spaces.  
\end{abstract}
\section{Introduction}
 \setcounter{equation}{0}
Let $0<\alpha<\N$. A fractional integral operator $I_\alpha$ is defined by
\bel{i_alpha}
\Big(I_\alpha f\Big)(x)~=~\int_{\R^\N} f(y)\left({1\over |x-y|}\right)^{\N-\alpha} dy.
\eeq
In 1928, Hardy and Littlewood \cite{Hardy-Littlewood} first established a weighted norm inequality for $I_\alpha$ in  one dimensional space, by considering the {\it weights} to be suitable power functions. This result has been extended to higher dimensions by Stein and Weiss \cite{Stein-Weiss}  and now bears the name of Stein-Weiss inequality.  
\v

{\bf Theorem A: Stein and Weiss (1958)}  ~{\it Let $\omega(x)=|x|^{-\gamma}, \sigma(x)=|x|^\delta,~\gamma, \delta\in\R$. We have
\bel{norm ineq}
\left\|\omega I_\alpha f\right\|_{\L^q\left(\R^\N\right)}~\leq~\C_{p~q~\alpha~\gamma~\delta~\N}~\left\| f\sigma\right\|_{\L^p\left(\R^\N\right)}
\eeq
for  $1<p\leq q<\infty$, if 
\bel{constraints}
\gamma~<~{\N\over q},\qquad \delta~<~\N\left({p-1\over p}\right),\qquad \gamma+\delta~\ge~0
\eeq
and
\bel{formula}
{\alpha\over \N}~=~{1\over p}-{1\over q}+{\gamma+\delta\over \N}.
\eeq}

$\diamond$ {\small Throughout, we regard $\C$ as a generic constant depending on its subindices}.

In the case of $\gamma=\delta=0$, {\bf Theorem A}  was proved in $\R^\N$ by Sobolev \cite{Sobolev}. This is known  today as Hardy-Littlewood-Sobolev inequality.

The weighted norm inequalities of fractional integrals have been extensively studied,  i.e:  by Muckenhoupt and Wheeden \cite{Muckenhoupt-Wheeden},  Coifman and Fefferman \cite{Coifman-C.Fefferman}, Fefferman and Muckenhoupt \cite{C.Fefferman-Muckenhoupt},      P\'{e}rez \cite{Perez} and Sawyer and Wheeden \cite{Sawyer-Wheeden}. 

Let $Q$ denote a cube in $\R^\N$. It is well known that  the norm inequality  (\ref{norm ineq}) implies 
\bel{A-characteristic}
\sup_{Q\subset\R^\N}~|Q|^{{\alpha\over \N}-\left({1\over p}-{1\over q}\right)}\left\{{1\over |Q|}\int_Q \omega^{q}(x)dx\right\}^{1\over q}\left\{{1\over|Q|}\int_Q \left({1\over \sigma}\right)^{p\over p-1}(x)dx\right\}^{p-1\over p}~<~\infty.
\eeq
The supremum (\ref{A-characteristic}) is  called the Muckenhoupt characteristic,  as was first introduced by Muckenhoupt for which $\omega^q$ and $\sigma^{-{p\over p-1}}$  are nonnegative and locally integrable functions. 

By taking into account  $\omega(x)=|x|^{-\gamma},\sigma(x)=|x|^\delta,~ \gamma, \delta\in\R$, we find that (\ref{A-characteristic}) implies the constraints  in (\ref{constraints})-(\ref{formula}). Hence,  (\ref{norm ineq}),   (\ref{constraints})-(\ref{formula}) and  (\ref{A-characteristic})  are equivalent conditions.

Consider $\R^\N$ as a product space, by writing 
$\R^\N=\R^{\N_1}\times\R^{\N_2}\times\cdots\times\R^{\N_n},~n\ge2$.
Let 
\bel{alpha}
0<\alpha_i<\N_i,\qquad i=1,2,\ldots,n\qquad\hbox{and}\qquad \alpha~=~\alpha_1+\alpha_2+\cdots+\alpha_n.
\eeq
In this paper,  we give an extension of {\bf Theorem A}  on  product spaces by studying so-called the  {\it strong fractional integral operator} $\I_\alpha$  defined by
\bel{I_alpha f}
\Big(\I_\alpha f\Big)(x)~=~\int_{\R^\N} f(y)\prod_{i=1}^n \left({1\over |x_i-y_i|}\right)^{\N_i-\alpha_i} dy,
\eeq
whose kernel  has singularity appeared on every coordinate subspace.

Study of certain operators  that  commute with a multi-parameter family of dilations,  dates back to the time of Jessen, Marcinkiewicz and Zygmund. During the past several decades, a number of pioneering  results have been accomplished, for example, 
by Robert Fefferman \cite{R.Fefferman}-\cite{R.Fefferman''}, Chang and Fefferman \cite{Chang-Fefferman},
Cordoba and Fefferman \cite{Cordoba-Fefferman},     Fefferman and Stein \cite{R.Fefferman-Stein},    M\"{u}ller, Ricci and Stein \cite{M.R.S},
Journ\'{e} \cite{Journe'} and Pipher \cite{Pipher}. The area remains largely open for  fractional integration.

\section{Statement of main result}
\setcounter{equation}{0}
{\bf Theorem A*:} {\it Let $\omega(x)=|x|^{-\gamma},\sigma(x)=|x|^\delta,~\gamma,\delta\in\R$. For $1<p\leq q<\infty$, the following  conditions  are equivalent: 

{\bf 1.} Let  $\Q\doteq\Q_1\times\Q_2\times\cdots\times\Q_n~\subset~\R^{\N_1}\times\R^{\N_2}\times\cdots\times\R^{\N_n}=\R^\N$ where $\Q_i$ denotes a cube in $\R^{\N_i}$, for every $i=1,2,\ldots,n$. 
\bel{A-Characteristic}
\sup_{\Q\subset\R^\N}~\prod_{i=1}^n |\Q_i|^{{\alpha_i\over \N_i}-\left({1\over p}-{1\over q}\right)}\left\{{1\over|\Q|}\int_{\Q} \omega^{q}\left(x\right)dx\right\}^{1\over q}\left\{{1\over|\Q|}\int_{\Q} \left({1\over \sigma}\right)^{p\over p-1}\left( x\right)dx\right\}^{p-1\over p}~<~\infty.
\eeq

{\bf 2.} 
\bel{local integrability}
\gamma~<~{\N\over q},\qquad \delta~<~\N\left({p-1\over p}\right),
\qquad \gamma+\delta~\ge~0
\eeq
and
\bel{Formula}
{\alpha\over\N}~=~{1\over p}-{1\over q}+{\gamma+\delta\over \N}.
\eeq 
For $\gamma\ge0,\delta\leq0$, 
\bel{Constraint Case One} 
\alpha_i-{\N_i\over p}~<~\delta,\qquad i=1,2,\ldots,n.
\eeq
For $\gamma\leq0,\delta\ge0$, 
\bel{Constraint Case Two}
 \alpha_i-\N_i\left({q-1\over q}\right)~<~\gamma,\qquad i=1,2,\ldots,n.
\eeq
For $\gamma>0,\delta>0$, 
\bel{Constraint Case Three}
\begin{array}{cc}\ds
\sum_{i\in\U} \alpha_i-{\N_i\over p}~<~\delta,\qquad
\U~=~\left\{i\in\{1,2,\ldots,n\}~\colon~ \alpha_i- {\N_i\over p}~\ge~0\right\},
\\\\ \ds
\sum_{i\in\V} \alpha_i-\left({q-1\over q}\right)\N_i~<~\gamma,\qquad \V~=~\left\{i\in\{1,2,\ldots,n\}~\colon~ \alpha_i-\N_i\left({q-1\over q}\right)~\ge~0\right\}.
\end{array}
\eeq

{\bf 3.} Let $\I_\alpha$ to be defined in (\ref{alpha})-(\ref{I_alpha f}). We have
\bel{Norm Ineq}
\left\|\omega\I_\alpha f\right\|_{\L^q\left(\R^\N\right)}~\leq~\C_{p~q~\alpha~\gamma~\delta~n~\N}~\left\|f\sigma\right\|_{\L^p\left(\R^\N\right)}.
\eeq}
\begin{remark}
In the $2$-parameter setting ($n=2$), {\bf Theorem A*} is first proved in the joint work by Sawyer and Wang \cite{Sawyer-Wang}. For $\gamma\ge0, \delta\leq0$ or $\gamma\leq0, \delta\ge0$, the "sandwiching" idea introduced in \cite{Sawyer-Wang} applies to the general multi-parameter situation. However, the difficult case occurs when $\gamma>0, \delta>0$, whereas the method used in  \cite{Sawyer-Wang} relies on solving a system of algebraic equations, which is no longer solvable for $n>2$.
\end{remark}
{\bf Sketch of Proof:} In section 3, we introduce a new framework, where the product space 
is decomposed into an infinitely many of dyadic cones. 
Every partial sum operator defined on a dyadic cone is essentially an one-parameter fractional integral operator, satisfying the desired regularity. 

In section 4, by taking into account $\omega(x)=|x|^{-\gamma},\sigma(x)=|x|^\delta,~ \gamma, \delta\in\R$, we prove that the Muckenhoupt characteristic  (\ref{A-Characteristic}) implies the constraints in (\ref{local integrability})-(\ref{Constraint Case Three}).

In section 5, by using  (\ref{local integrability})-(\ref{Constraint Case Three}), we show that 
\bel{A-Characteristic r-bump}
\prod_{i=1}^n |\Q_i|^{{\alpha_i\over \N_i}-\left({1\over p}-{1\over q}\right)}\left\{{1\over|\Q|}\int_{\Q} \omega^{qr}\left(x\right)dx\right\}^{1\over qr}\left\{{1\over|\Q|}\int_{\Q} \left({1\over \sigma}\right)^{pr\over p-1}\left( x\right)dx\right\}^{p-1\over pr},\qquad r>1
\eeq
decays exponentially,  as  the eccentricity of $\Q$ getting large, for
$\alpha_i>\N_i\left({1\over p}-{1\over q}\right),~i=1,2,\ldots,n$.

On the other hand,  we handle  the case 
$\alpha_i=\N_i\left({1\over p}-{1\over q}\right),~i=1,2,\ldots,n$ in section 6. 
  
We prove {\bf Theorem A*} in the last section, by decomposing $\I_\alpha$ so that the resulting estimates can be reduced to either of the above two cases.

For dealing with such convolution operators with positive kernels, it is suffice to assume $f\ge0$ in the rest of the paper.

\section{Cone decomposition on product spaces}
\setcounter{equation}{0}
Let $\t$  denote an $n$-tuple $(2^{-t_1},2^{-t_2},\ldots,2^{-t_n})$ where  $t_i,~i=1,2,\ldots,n$  are nonnegative integers. We require
$t_\nu\doteq\min\{t_i\colon~i=1,2,\ldots,n\}=0$.

Define
\bel{Partial}
\Big(\Delta_\t \I_\alpha f\Big)(x)~\doteq~\int_{\Gamma_\t(x)} f(y)\prod_{i=1}^n \left({1\over |x_i-y_i|}\right)^{\N_i-\alpha_i}dy
\eeq
where
\bel{Cone}
\begin{array}{lr}\ds
\Gamma_\t(x)~\doteq~\bigotimes_{i=1}^n\left\{y_i\in\R^{\N_i}\colon~2^{-t_i}\leq {|x_i-y_i|\over|x_\nu-y_\nu|}< 2^{-t_i+1}\right\}.
\end{array}
\eeq
Observe that $\Gamma_\t(x)$ in (\ref{Cone}) is a dyadic cone with vertex on $x$ whose eccentricity depends on $\t$. In particular, we write 
\bel{Gamma_o}
\Gamma_o(x)~\doteq~\Gamma_\t(x),\qquad t_1=t_2=\cdots=t_n=0.
\eeq
\begin{figure}[h]
\centering
\includegraphics[scale=0.52]{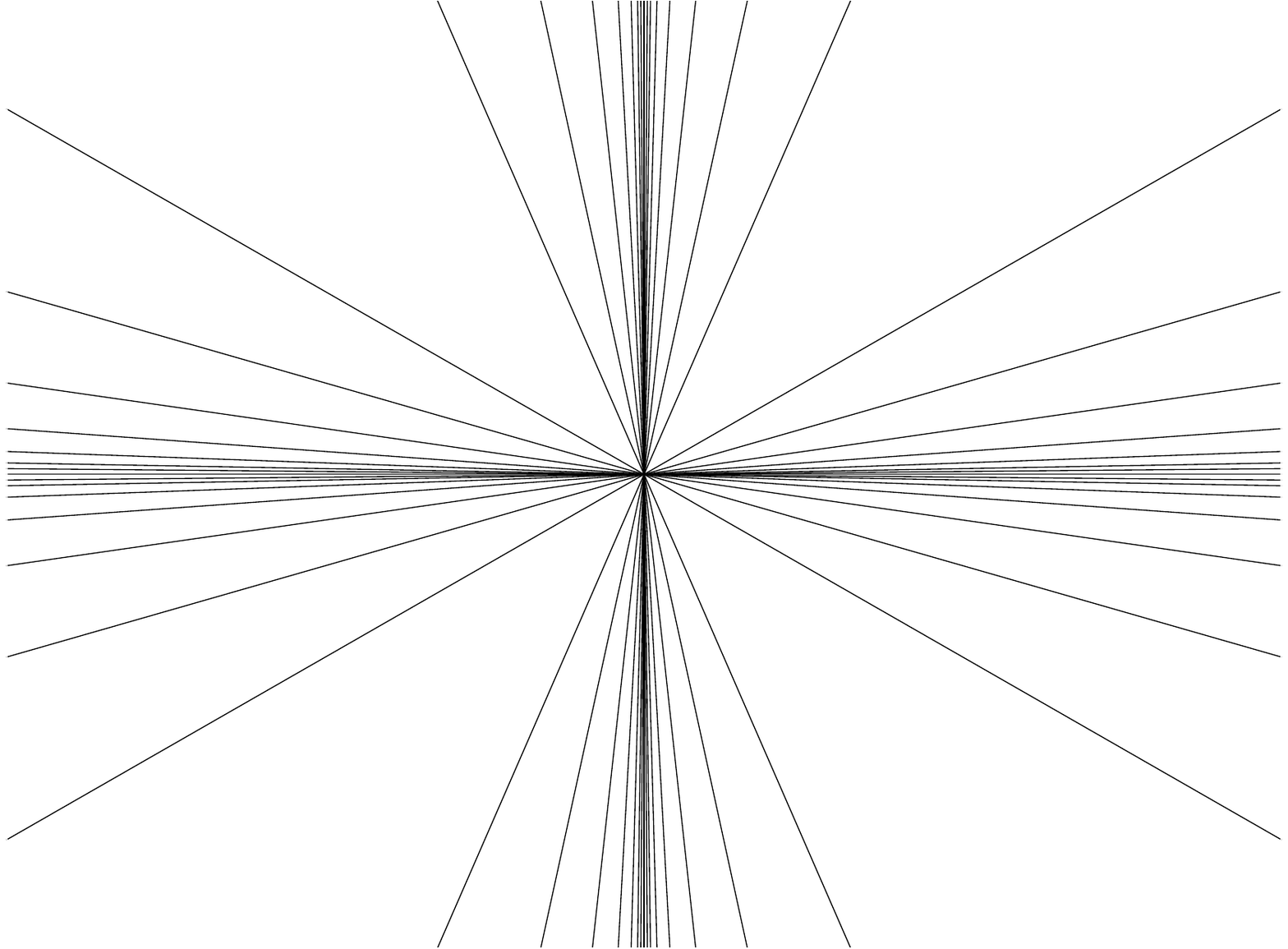}
\caption{dyadic cones in a $2$-parameter setting.}
\end{figure}

Denote an $n$-parameter dilation
\bel{dila}
\t x~=~\left( 2^{-t_1}x_1,2^{-t_2}x_2,\dots,2^{-t_n}x_n\right).
\eeq
Let $\Q^\t$ be a dilated of $\Q$ such that
$|\Q^\t_i|^{1\over \N_i}=2^{-t_i}|\Q_i|^{1\over \N_i}, i=1,2,\ldots,n$.
We have
\bel{A-Characteristic Dila}
\begin{array}{lr}\ds
~~~~~~~\prod_{i=1}^n |\Q_i|^{{\alpha_i\over \N_i}-\left({1\over p}-{1\over q}\right)}\left\{{1\over|\Q|}\int_{\Q} \omega^{qr}\left(\t x\right)dx\right\}^{1\over qr}\left\{{1\over|\Q|}\int_{\Q} \left({1\over \sigma}\right)^{pr\over p-1}\left(\t x\right)dx\right\}^{p-1\over pr}
\\\\ \ds
~=~\prod_{i=1}^n |\Q_i|^{{\alpha_i\over \N_i}-\left({1\over p}-{1\over q}\right)}\left\{{1\over|\Q^\t|}\int_{\Q^\t} \omega^{qr}\left( x\right)dx\right\}^{1\over qr}\left\{{1\over|\Q^\t|}\int_{\Q^\t} \left({1\over \sigma}\right)^{pr\over p-1}\left( x\right)dx\right\}^{p-1\over pr}
\\\\ \ds
~=~\prod_{i=1}^n2^{t_i\left(\alpha_i-{\N_i\over p}+{\N_i\over q}\right)}\prod_{i=1}^n  |\Q_i^\t|^{{\alpha_i\over \N_i}-\left({1\over p}-{1\over q}\right)}\left\{{1\over|\Q^\t|}\int_{\Q^\t} \omega^{qr}\left( x\right)dx\right\}^{1\over qr}\left\{{1\over|\Q^\t|}\int_{\Q^\t} \left({1\over \sigma}\right)^{pr\over p-1}\left( x\right)dx\right\}^{p-1\over pr}
\end{array}
\eeq
for every $\Q\subset\R^\N$.

Given $\t$, consider 
\bel{Q ratio}
{^\t}\Q~\subset~\R^\N~\colon~{|\Q_i|^{1\over \N_i}/ |\Q_\nu|^{1\over\N_\nu}}~=~2^{-t_i},~~~  i=1,2,\ldots,n.
\eeq
For $r\ge1$, we define
\bel{sup A_pqr^alpha t}
\begin{array}{lr}\ds
\A_{pqr}^\alpha\left(\t~\colon\omega,\sigma\right)~=~
\sup_{{^\t}\Q} ~\prod_{i=1}^n |\Q_i|^{{\alpha_i\over \N_i}-\left({1\over p}-{1\over q}\right)}\left\{{1\over|\Q|}\int_\Q \omega^{qr}\left(x\right)dx\right\}^{1\over qr}\left\{{1\over|\Q|}\int_\Q \left({1\over \sigma}\right)^{pr\over p-1}\left( x\right)dx\right\}^{p-1\over pr}.
\end{array}
\eeq
Suppose that $\Q$ satisfies
$\left|\Q_1\right|^{1\over \N_1}=\left|\Q_2\right|^{1\over\N_2}=\cdots=\left|\Q_n\right|^{1\over \N_n}$.
We have $\Q^\t={^\t}\Q$ and
\bel{A-Characteristic Bound}
\begin{array}{lr}\ds
~~~~~~~ |\Q|^{{\alpha\over \N}-\left({1\over p}-{1\over q}\right)}\left\{{1\over|\Q|}\int_{\Q} \omega^{qr}\left(\t x\right)dx\right\}^{1\over qr}\left\{{1\over|\Q|}\int_{\Q} \left({1\over \sigma}\right)^{pr\over p-1}\left(\t x\right)dx\right\}^{p-1\over pr}
\\\\ \ds
~=~\prod_{i=1}^n2^{t_i\left(\alpha_i-{\N_i\over p}+{\N_i\over q}\right)} \prod_{i=1}^n |\Q_i^\t|^{{\alpha_i\over \N_i}-\left({1\over p}-{1\over q}\right)}\left\{{1\over|\Q^\t|}\int_{\Q^\t} \omega^{qr}\left( x\right)dx\right\}^{1\over qr}\left\{{1\over|\Q^\t|}\int_{\Q^\t} \left({1\over \sigma}\right)^{pr\over p-1}\left( x\right)dx\right\}^{p-1\over pr}~~\hbox{\small{by (\ref{A-Characteristic Dila})}}
\\\\ \ds
~\leq~\prod_{i=1}^n2^{t_i\left(\alpha_i-{\N_i\over p}+{\N_i\over q}\right)} \A_{pqr}^\alpha\left(\t~\colon\omega,\sigma\right)\qquad\hbox{\small{by (\ref{Q ratio})-(\ref{sup A_pqr^alpha t}).}}
\end{array}
\eeq
Now, recall  Sawyer-Wheeden theorem for  one-parameter fractional integral operators in weighted norms, stated as Theorem 1 in \cite{Sawyer-Wheeden}:
\bel{One-para Norm Ineq}
\begin{array}{lr}\ds
\left\{\int_{\R^\N}\left\{\int_{\R^\N}f( y)\left({1\over |x-y|}\right)^{\N-\alpha} dy\right\}^{q}\omega^q(x)dx\right\}^{1\over q} 
~\leq~\C_{p~q~r~\alpha~\N}~ A_{pqr}^\alpha(\omega, \sigma)~\left\{\int_{\R^\N} \Big(f\sigma\Big)^p(x)dx\right\}^{1\over p}
\end{array}
\eeq
for $1<p\leq q<\infty$, if 
\bel{A-Chara r-bump}
\begin{array}{lr}\ds
A_{pqr}^\alpha(\omega, \sigma)~\doteq~
\\ \ds
\sup_{\Q~\colon~\left|\Q_1\right|^{1\over \N_1}~=~\cdots~=~\left|\Q_n\right|^{1\over \N_n}} |\Q|^{{\alpha\over \N}-\left({1\over p}-{1\over q}\right)}\left\{{1\over |\Q|}\int_\Q \omega^{qr}(x)dx\right\}^{1\over qr}\left\{{1\over|\Q|}\int_\Q \left({1\over \sigma}\right)^{pr\over p-1}(x)dx\right\}^{p-1\over pr}~<~\infty
\\ \ds~~~~~~~~~~~~~~~~~~~~~~~~~~~~~~~~~~~~~~~~~~~~~~~~~~~~~~~~~~~~~~~~~~~~~~~~~~~~~~~~~~~~~~~~~~~~~~~~~~~~~~~~~~~~~~~~~~~~~~~~~~~~~~~~~r>1.
\end{array}
\eeq

\begin{remark}The constant $\C_{p~q~r~\alpha~\N}~ A_{pqr}^\alpha(\omega, \sigma)$ in (\ref{One-para Norm Ineq}) is not written explicitly in the statement of Theorem 1 by Sawyer and Wheeden \cite{Sawyer-Wheeden}. But it can be computed directly by carrying out the proof given in section 2 of \cite{Sawyer-Wheeden}.\end{remark}
By applying (\ref{One-para Norm Ineq})-(\ref{A-Chara r-bump})  and  using the estimate in (\ref{A-Characteristic Bound}), we have
\bel{Regularity est}
\begin{array}{lr}\ds
\left\{\int_{\R^\N}\left\{\int_{\R^\N}f(\t y)\left({1\over |x-y|}\right)^{\N-\alpha} dy\right\}^{q}\omega^q(\t x)dx\right\}^{1\over q} 
\\\\ \ds
~\leq~\C_{p~q~r~\alpha~\N}~\prod_{i=1}^n2^{t_i\left(\alpha_i-{\N_i\over p}+{\N_i\over q}\right)} \A_{pqr}^\alpha\left(\t~\colon\omega,\sigma\right)\left\{\int_{\R^\N} \Big(f\sigma\Big)^p(\t x)dx\right\}^{1\over p}
\end{array}
\eeq
for $1<p\leq q<\infty$ and every $\t$.

Recall from (\ref{Partial})-(\ref{Cone}). By changing dilations $x\mt \t x, y\mt \t y$, we have
\bel{Dilation Est}
\begin{array}{lr}\ds
\left\{\int_{\R^\N} \Big(\Delta_\t\I_{\alpha}f\Big)^q(x)\omega^q(x)dx\right\}^{1\over q}
\\\\ \ds
~=~\left\{\int_{\R^\N}\left\{ \int_{\Gamma_\t(x)}f(y)\prod_{i=1}^n\left({1\over |x_i-y_i|}\right)^{\N_i-\alpha_i}dy\right\}^q\omega^q(x)dx\right\}^{1\over q}
 \\\\ \ds
 ~=~\left\{\int_{\R^\N}\left\{\int_{\Gamma_o(x)}f\left(\t y\right)\left\{\prod_{i=1}^n2^{-t_i\N_i}\left({1\over 2^{-t_i}|x_i-y_i|}\right)^{\N_i-\alpha_i} \right\}dy\right\}^q\omega^q\left(\t x\right)\prod_{i=1}^n 2^{-t_i\N_i} dx\right\}^{1\over q}
 \\\\ \ds
 ~\leq~\C_{\alpha~n~\N}~\prod_{i=1}^n2^{-t_i\left(\alpha_i+{\N_i\over q}\right)}\left\{\int_{\R^\N}\left\{ \int_{\R^\N}f\left(\t y\right)\left({1\over |x-y|}\right)^{\N-\alpha}dy\right\}^q\omega^q\left(\t x\right)dx\right\}^{1\over q}
 \\\\ \ds
 ~\leq~\C_{p~q~r~\alpha~n~\N}~\prod_{i=1}^n 2^{-t_i\left(\alpha_i+{\N_i\over q}\right)}2^{t_i\left(\alpha_i-{\N_i\over p}+{\N_i\over q}\right)}\A_{pqr}^\alpha\left(\t~\colon\omega,\sigma\right)\left\{ \int_{\R^\N} \Big(f\sigma\Big)^p\left(\t x\right)dx\right\}^{1\over p}\qquad \hbox{\small{by (\ref{Regularity est})}}
  \\\\ \ds
   ~=~\C_{p~q~r~\alpha~n~\N}~\A_{pqr}^\alpha\left(\t~\colon\omega,\sigma\right)\prod_{i=1}^n 2^{-t_i\left(\alpha_i+{\N_i\over q}\right)}2^{t_i\left(\alpha_i-{\N_i\over p}+{\N_i\over q}\right)}\left\{ \int_{\R^\N} \Big(f\sigma\Big)^p\left(x\right)\prod_{i=1}^n 2^{t_i\N_i}dx\right\}^{1\over p}
 \\\\ \ds
 ~=~\C_{p~q~r~\alpha~n~\N}~\A_{pqr}^\alpha\left(\t~\colon\omega,\sigma\right)\left\{ \int_{\R^\N} \Big(f\sigma\Big)^p\left(x\right)dx\right\}^{1\over p}.
\end{array}
\eeq
Observe that $\Delta_\t\I_\alpha$ is essentially an one-parameter fractional integral operator,   satisfying 
\bel{Partial Result}
\left\|\Big(\Delta_\t \I_\alpha f\Big)\omega\right\|_{\L^q\left(\R^\N\right)}~\leq~\C_{p~q~r~\alpha~n~\N}~\A_{pqr}^\alpha\left(\t~\colon\omega,\sigma\right)\left\| f\sigma\right\|_{\L^p\left(\R^\N\right)}
\eeq
for $1<p\leq q<\infty$. 

By applying Minkowski inequality,   provided that
\bel{Summability}
\sum_\t \A_{pqr}^\alpha\left(\t~\colon\omega,\sigma\right)~<~\infty,
\eeq
the norm inequality holds in (\ref{Norm Ineq}).

\section{Necessary constraints}
\setcounter{equation}{0}
First, it is well known that the norm inequality (\ref{Norm Ineq}) implies 
\bel{A-Characteristic sup}
\begin{array}{lr}\ds
\A_{pq}^\alpha(\omega,\sigma)~\doteq~\sup_{\Q\subset\R^\N}~\prod_{i=1}^n |\Q_i|^{{\alpha_i\over \N_i}-\left({1\over p}-{1\over q}\right)}\left\{{1\over|\Q|}\int_\Q \omega^q(x)dx\right\}^{1\over q}\left\{{1\over|\Q|}\int_\Q \left({1\over \sigma}\right)^{p\over p-1}(x)dx\right\}^{p-1\over p}~<~\infty.
\end{array}
\eeq
Let $\omega(x)=|x|^{-\gamma},\sigma(x)=|x|^\delta, \gamma,\delta\in\R$. We aim to show the Muckenhoupt characteristic (\ref{A-Characteristic sup}) implying the constraints in  (\ref{local integrability})-(\ref{Constraint Case Three}).

Let $\Q^\lambda$ denote a dilated variant of $\Q$ for $\lambda>0$,   such that $\Q^\lambda\doteq\Q^\lambda_1\times\Q^\lambda_2\times\cdots\times\Q^\lambda_n$ and $|\Q_i^\lambda|^{1\over \N_i}=\lambda |\Q_i|^{1\over \N_i},~ i=1,2,\ldots,n$. Suppose $\omega(x)=|x|^{-\gamma},\sigma(x)=|x|^\delta, \gamma,\delta\in\R$. From (\ref{A-Characteristic sup})-(\ref{A-Characteristic*}), we have
\bel{Equi Chara}
\begin{array}{lr}\ds
~~~~~~~\prod_{i=1}^n |\Q_i|^{{\alpha_i\over \N_i}-\left({1\over p}-{1\over q}\right)}\left\{{1\over |\Q|}\int_\Q \left({1\over |x|}\right)^{\gamma q}dx\right\}^{1\over q}\left\{{1\over|\Q|}\int_\Q \left({1\over |x|}\right)^{\delta p\over p-1}dx\right\}^{p-1\over p}
\\\\ \ds
~=~\lambda^{\gamma+\delta-\alpha+\N\left({1\over p}-{1\over q}\right)}\prod_{i=1}^n |\Q_i^\lambda|^{{\alpha_i\over \N_i}-\left({1\over p}-{1\over q}\right)}\left\{{1\over |\Q^\lambda|}\int_{\Q^\lambda} \left({1\over |x|}\right)^{\gamma q}dx\right\}^{1\over q}\left\{{1\over|\Q^\lambda|}\int_{\Q^\lambda} \left({1\over |x|}\right)^{\delta p\over p-1}dx\right\}^{p-1\over p}
\\\\ \ds
~\leq~\lambda^{\gamma+\delta-\alpha+\N\left({1\over p}-{1\over q}\right)}\A_{pq}^\alpha\left(|x|^{-\gamma},|x|^\delta\right)~<~\infty.
\end{array}
\eeq
Consider $|\Q_1|^{1\over\N_1}=|\Q_2|^{1\over\N_2}=\cdots=|\Q_n|^{1\over\N_n}=1$. The first line of (\ref{Equi Chara}) is bounded from below. Suppose $\gamma+\delta-\alpha+\N\left({1\over p}-{1\over q}\right)\neq0$.  By either taking $\lambda\mt0$ or $\lambda\mt\infty$, the last line of (\ref{Equi Chara})  is  vanished. Hence that  we must have
$\gamma+\delta-\alpha+\N\left({1\over p}-{1\over q}\right)=0$
which is (\ref{Formula}).

We write $x=(x_i,x_i^\dagger)\in\R^{\N_i}\times\R^{\N-\N_i}, ~i=1,2,\ldots,n$ and $\Q_i^\dagger=\bigotimes_{j\neq i}\Q_j$.
Let $\Q_i$  shrink to some $x_i\in\Q_i$ and  $|\Q_j|^{1\over \N_j}=1,~j\neq i$ in (\ref{A-Characteristic sup}). 
Suppose $x_i\neq0$ in $\R^{\N_i}$. By applying Lebesgue Differentiation Theorem, we have
\bel{A-Characteristic-i}
\begin{array}{lr}\ds
\left\{\lim_{|\Q_i|\mt 0} |\Q_i|^{{\alpha_i\over \N_i}-\left({1\over p}-{1\over q}\right)}\right\}\left\{
{1\over |\Q_i^\dagger|}\int_{\Q_i^\dagger} \Bigg({1\over |x_i+x_i^\dagger|}\Bigg)^{\gamma q} dx_i^\dagger\right\}^{1\over q}\left\{{1\over|\Q_i^\dagger|}\int_{\Q_i^\dagger}  \Bigg({1\over |x_i+x_i^\dagger|}\Bigg)^{\delta p\over p-1}dx_i^\dagger\right\}^{p-1\over p}
\\\\ \ds
~\leq~\A_{pq}^\alpha\left(|x|^{-\gamma},|x|^\delta\right),~~~~~~~~~~i~=~1,2,\ldots,n.
\end{array}
\eeq
Note that $|\Q_i^\dagger|=1$ in (\ref{A-Characteristic-i}). The boundedness of $\A_{pq}^\alpha\left(|x|^{-\gamma},|x|^\delta\right)$ requires
\bel{subbalance}
{\alpha_i\over\N_i}~\ge~{1\over p}-{1\over q},\qquad i~=~1,2,\ldots,n.
\eeq
By putting together (\ref{subbalance}) and (\ref{Formula}), we find $\gamma+\delta\ge0$.  On the other and, it is essential to require
$\gamma q<\N$ and $\delta\left({p\over p-1}\right)<\N$
for the local integrability of $|x|^{-\gamma q}$ and $|x|^{-\delta\left({p\over p-1}\right)}$ respectively.
These are the constraints in (\ref{local integrability}).

In the remaining section, we assume $\Q$ centered on the origin of $\R^\N$.
 
Let $\S$ to be a proper subset of $\{1,2,\ldots,n\}$. We define the truncated cube  $\Q_i^\ve=\Q_i\cap\{|x_i|\ge\ve\}$ for  $\ve>0$ and every $i\in\S$. 
Denote $\Q^\ve\doteq\bigotimes_{i\in\S}\Q_i^\ve\times \bigotimes_{i\in\S^c}\Q_i$ and 
$\Q_{\S}=\bigotimes_{i\in\S}\Q_i,~ \Q_{\S^c}=\bigotimes_{i\in\S^c}\Q_i$. Moreover,  we write $x=\left(x_\S,x_{\S^c}\right)\in\R^{\N_\S}\times\R^{\N-\N_\S}$ for which $\N_\S=\sum_{i\in\S}\N_i$

Suppose that there exists at least one $i\in\S^c$ such that 
$\alpha_i-\N_i\left({1\over p}-{1\over q}\right)>0$. 
Let $0<\lambda<1$. Consider $|\Q_i|^{1\over \N_i}=1$ for $i\in\S$ and $|\Q_i|^{1\over \N_i}=\lambda$ for   $i\in\S^c$. We have
\bel{A_pq Decay Est}
\begin{array}{lr}\ds
\prod_{i=1}^n |\Q_i|^{{\alpha_i\over \N_i}-\left({1\over p}-{1\over q}\right)}\left\{{1\over|\Q|}\int_{\Q} \left({1\over|x|}\right)^{\gamma q}dx\right\}^{1\over q}\left\{{1\over|\Q|}\int_{\Q} \left({1\over |x|}\right)^{\delta p\over p-1}dx\right\}^{p-1\over p}
\\\\ \ds
~=~\lim_{\ve\mt0}\prod_{i=1}^n |\Q_i|^{{\alpha_i\over \N_i}-\left({1\over p}-{1\over q}\right)}\left\{{1\over|\Q|}\int_{\Q^\ve} \left({1\over|x|}\right)^{\gamma q}dx\right\}^{1\over q}\left\{{1\over|\Q|}\int_{\Q^\ve} \left({1\over |x|}\right)^{\delta p\over p-1}dx\right\}^{p-1\over p}
\\\\ \ds
~=~\lim_{\ve\mt0}~\lambda^{\sum_{i\in\S^c} \alpha_i-\N_i\left({1\over p}-{1\over q}\right)} \left\{ {1\over|\Q|}\int_{\Q^\ve} \left({1\over|x|}\right)^{\gamma q}dx\right\}^{1\over q}\left\{{1\over|\Q|}\int_{\Q^\ve} \left({1\over |x|}\right)^{\delta p\over p-1}dx\right\}^{p-1\over p}

\\\\ \ds
~=~\lim_{\ve\mt0}~ 0\times
\left\{ \idotsint_{\bigotimes_{i\in\S} \Q_i^\ve}\left({1\over\sum_{i\in\S}|x_i|^2}\right)^{\gamma q\over 2}\prod_{i\in\S}dx_i\right\}^{1\over q}
\left\{  \idotsint_{\bigotimes_{i\in\S}\Q_i^\ve} \left({1\over \sum_{i\in\S}|x_i|^2}\right)^{{1\over 2}{\delta p\over p-1}}\prod_{i\in\S} dx_i\right\}^{p-1\over p}
\\ \ds~~~~~~~~~~~~~~~~~~~~~~~~~~~~~~~~~~~~~~~~~~~~~~~~~~~~~~~~~~~~~~~~~~~~~~~~~~~~~~~~~ \hbox{\small{by Lebesgue differentiation theorem at $\lambda=0$}}
\\\\ \ds
~=~\lim_{\ve\mt0} 0~=~0.
\end{array}
\eeq
Suppose $\alpha_i-\N_i\left({1\over p}-{1\over q}\right)=0$ for every $i\in\S^c$. Let    $\Q_i$ shrink to the origin of $\R^{\N_i}$ for  every $i\in\S^c$  in (\ref{A-Characteristic sup}). 
By applying Lebesgue differentiation theorem, we have
\bel{A-subspace S}
\begin{array}{lr}\ds
\A_{pq}^\alpha\left(|x|^{-\gamma},|x|^\delta\right)~\ge~
\prod_{i=1}^n |\Q_i|^{{\alpha_i\over \N_i}-\left({1\over p}-{1\over q}\right)}\left\{{1\over |\Q|}\int_\Q \left({1\over |x|}\right)^{\gamma q}dx\right\}^{1\over q}\left\{{1\over|\Q|}\int_\Q \left({1\over |x|}\right)^{\delta p\over p-1}dx\right\}^{p-1\over p}
\\\\ \ds~~~~~~~~~~~~~~~~~~~~~~~~
~=~\prod_{i\in\S} |\Q_i|^{{\alpha_i\over \N_i}-\left({1\over p}-{1\over q}\right)}\left\{{1\over|\Q_\S|}\int_{\Q_\S} \left({1\over|x_\S|}\right)^{\gamma q}dx_\S\right\}^{1\over q}\left\{{1\over|\Q_\S|}\int_{\Q_\S} \left({1\over |x_\S|}\right)^{\delta p\over p-1}dx_\S\right\}^{p-1\over p}
\end{array}
\eeq
where $\gamma q<\N_\S$ and $\delta\left({p\over p-1}\right)<\N_\S$ become necessities.

{\bf Case One:} Consider $\gamma\ge0, \delta\leq0$.  
Let  $|\Q_i|^{1\over \N_i}=1$ for $i\in\{1,2,\ldots,n\}$ and  $|\Q_j|^{1\over \N_j}=\lambda$ for all $j\neq i$. 
Suppose  $\alpha_j-\N_j\left({1\over p}-{1\over q}\right)=0$ for every $j\neq i$. We have
\bel{Characteristic Est1 Case1}
\begin{array}{lr}\ds
\prod_{i=1}^n |\Q_i|^{{\alpha_i\over \N_i}-\left({1\over p}-{1\over q}\right)}\left\{{1\over|\Q|}\int_{\Q} \left({1\over|x|}\right)^{\gamma q}dx\right\}^{1\over q}\left\{{1\over|\Q|}\int_{\Q} \left({1\over |x|}\right)^{\delta p\over p-1}dx\right\}^{p-1\over p}
\\\\ \ds
~\ge~\C_{q~\gamma~n}~\left\{\int_{\Q_i}\left({1\over \lambda+|x_i|}\right)^{\gamma q}dx_i\right\}^{1\over q}
\left\{\int_{\Q_i}\left({1\over |x_i|}\right)^{\delta p\over p-1}dx_i\right\}^{p-1\over p}\qquad \hbox{\small{($\delta\leq0$)}}
\\\\ \ds
~\ge~\C_{p~q~\gamma~\delta~n~\N}~\left\{\int_{\lambda<|x_i|\leq1}\left({1\over \lambda+|x_i|}\right)^{\gamma q}dx_i\right\}^{1\over q}
\end{array}
\eeq
where
\bel{Necessary Case1 Est1}
\begin{array}{cc}\ds
\int_{\lambda<|x_i|\leq1}\left({1\over \lambda+|x_i|}\right)^{\gamma q}dx_i~\leq~\C_\N~\ln\left({1+\lambda\over 2\lambda}\right)~~~~\hbox{if}~~~~\gamma={\N_i\over q},
\\\\ \ds
\int_{\lambda<|x_i|\leq1}\left({1\over \lambda+|x_i|}\right)^{\gamma q}dx_i~\leq~\C_\N~{1\over \gamma q-\N_i}\left[ \left({1\over 2\lambda}\right)^{\gamma q-\N_i}-\left({1\over \lambda+1}\right)^{\gamma q-\N_i}\right]~~~~\hbox{if}~~~~\gamma~>~{\N_i\over q}.
\end{array}
\eeq
From (\ref{Characteristic Est1 Case1})-(\ref{Necessary Case1 Est1}), as $\lambda\mt0$, we need
 \bel{Necessity est1} 
\gamma~<~{\N_i\over q},\qquad i~=1,2,\ldots,n
\eeq
in order to satisfy the inequality in (\ref{Equi Chara}).

Suppose that there exists $j\neq i$ such that $\alpha_j-\N_j\left({1\over p}-{1\over q}\right)>0$. 
We have
\bel{Characteristic Est2 Case1}
\begin{array}{lr}\ds
\prod_{i=1}^n |\Q_i|^{{\alpha_i\over \N_i}-\left({1\over p}-{1\over q}\right)}\left\{{1\over|\Q|}\int_{\Q} \left({1\over|x|}\right)^{\gamma q}dx\right\}^{1\over q}\left\{{1\over|\Q|}\int_{\Q} \left({1\over |x|}\right)^{\delta p\over p-1}dx\right\}^{p-1\over p}
\\\\ \ds
~\ge~\C_{q~\gamma~n}\prod_{j\neq i}\lambda^{\alpha_j-\N_j\left({1\over p}-{1\over q}\right)}\left\{\int_{\Q_i}\left({1\over \lambda+|x_i|}\right)^{\gamma q}dx_i\right\}^{1\over q}
\left\{\int_{\Q_i}\left({1\over |x_i|}\right)^{\delta p\over p-1}dx_i\right\}^{p-1\over p}~~ \hbox{\small{($\delta\leq0$)}}
\\\\ \ds
~\ge~\C_{p~q~\gamma~\delta~n~\N}\prod_{j\neq i}\lambda^{\alpha_j-\N_j\left({1\over p}-{1\over q}\right)}\left\{\int_{0<|x_i|\leq\lambda}\left({1\over \lambda}\right)^{\gamma q}dx_i\right\}^{1\over q}
\\\\ \ds
~=~\C_{p~q~\gamma~\delta~n~\N}~\lambda^{{\N_i\over q}-\gamma+\sum_{j\neq i} \alpha_j-\N_j\left({1\over p}-{1\over q}\right)}.
\end{array}
\eeq
Recall the estimate in (\ref{A_pq Decay Est}) and take $\S=\{i\}$. We have (\ref{Characteristic Est2 Case1}) equal to zero at $\lambda=0$. Together with (\ref{Necessity est1}),  we find
\bel{Necessary Est1}
\gamma~<~{\N_i\over q}+\sum_{j\neq i} \alpha_j-\N_j\left({1\over p}-{1\over q}\right)
,\qquad i~=~1,2,\ldots,n.
\eeq
\v

{\bf Case Two:} Consider $\gamma\leq0, \delta\ge0$. 
Let $|\Q_i|^{1\over \N_i}=1$ for $i\in\{1,2,\ldots,n\}$ and  $|\Q_j|^{1\over \N_j}=\lambda$ for all  $j\neq i$.
Suppose  $\alpha_j-\N_j\left({1\over p}-{1\over q}\right)=0$ for every $j\neq i$. We have
\bel{Characteristic Est1 Case2}
\begin{array}{lr}\ds
\prod_{i=1}^n |\Q_i|^{{\alpha_i\over \N_i}-\left({1\over p}-{1\over q}\right)}\left\{{1\over|\Q|}\int_{\Q} \left({1\over|x|}\right)^{\gamma q}dx\right\}^{1\over q}\left\{{1\over|\Q|}\int_{\Q} \left({1\over |x|}\right)^{\delta p\over p-1}dx\right\}^{p-1\over p}
\\\\ \ds
~\ge~\C_{p~\delta~n}\left\{\int_{\Q_i}\left({1\over |x_i|}\right)^{\gamma q}dx_i\right\}^{1\over q}
\left\{\int_{\Q_i}\left({1\over \lambda+|x_i|}\right)^{\delta p\over p-1}dx_i\right\}^{p-1\over p}\qquad \hbox{\small{($\gamma\leq0$)}}
\\\\ \ds
~\ge~\C_{p~q~\gamma~\delta~n~\N}~\left\{\int_{\lambda<|x_i|\leq1}\left({1\over \lambda+|x_i|}\right)^{\delta p\over p-1}dx_i\right\}^{p-1\over p}
\end{array}
\eeq
where
\bel{Necessary Case2 Est1}
\begin{array}{cc}\ds
\int_{\lambda<|x_i|\leq1}\left({1\over \lambda+|x_i|}\right)^{\delta \left({p\over p-1}\right)}dx_i~\leq~\C_\N~\ln\left({1+\lambda\over2\lambda}\right)~~~~\hbox{if}~~~~\delta=\N_i\left({p-1\over p}\right),
\\\\ \ds
\int_{\lambda<|x_i|\leq1}\left({1\over \lambda+|x_i|}\right)^{\delta \left({p\over p-1}\right)}dx_i~\leq~\C_\N~{1\over  \delta \left({p\over p-1}\right)-\N_i} \left[\left({1\over 2\lambda}\right)^{\delta \left({p\over p-1}\right)-\N_i}-\left({1\over \lambda+1}\right)^{\delta \left({p\over p-1}\right)-\N_i}\right]
\\ \\ \ds
~~~~~~~~~~~~~~~~~~~~~~~~~~~~~~~~~~~~~~~~~~~~~~~~~~~~~~~~~~~~~~~~~~~~~~~~~~~~~~~~~~~~~~~~~~~~~~~~~~~~~~~\hbox{if}~~ \delta~>~\N_i\left({p-1\over p}\right).
\end{array}
\eeq
From (\ref{Characteristic Est1 Case2})-(\ref{Necessary Case2 Est1}), as $\lambda\mt0$, we need
 \bel{Necessity est2} 
\delta~<~\N_i\left({p-1\over p}\right),\qquad i~=1,2,\ldots,n
\eeq
in order to satisfy the inequality in (\ref{Equi Chara}).

Suppose that there exists  $j\neq i$ such that $\alpha_j-\N_j\left({1\over p}-{1\over q}\right)>0$. 
We have
\bel{Characteristic Est2 Case2}
\begin{array}{lr}\ds
\prod_{i=1}^n |\Q_i|^{{\alpha_i\over \N_i}-\left({1\over p}-{1\over q}\right)}\left\{{1\over|\Q|}\int_\Q \left({1\over|x|}\right)^{\gamma q}dx\right\}^{1\over q}\left\{{1\over|\Q|}\int_\Q \left({1\over |x|}\right)^{\delta p\over p-1}dx\right\}^{p-1\over p}
\\\\ \ds
~\ge~\C_{p~\delta~n}~\prod_{j\neq i}\lambda^{\alpha_j-\N_j\left({1\over p}-{1\over q}\right)}\left\{\int_{\Q_i}\left({1\over |x_i|}\right)^{\gamma q}dx_i\right\}^{1\over q}
\left\{\int_{\Q_i}\left({1\over \lambda+|x_i|}\right)^{\delta p\over p-1}dx_i\right\}^{p-1\over p}~~ \hbox{\small{($\gamma\leq0$)}}
\\\\ \ds
~\ge~\C_{p~q~\gamma~\delta~n~\N}\prod_{j\neq i}\lambda^{\alpha_j-\N_j\left({1\over p}-{1\over q}\right)}
\left\{\int_{0<|x_i|\leq\lambda}\left({1\over \lambda}\right)^{\delta p\over p-1}dx_i\right\}^{p-1\over p}
\\\\ \ds
~=~\C_{p~q~\gamma~\delta~n~\N}~\lambda^{\left({p-1\over p}\right)\N_i-\delta+\sum_{j\neq i} \alpha_j-\N_j\left({1\over p}-{1\over q}\right)}.
\end{array}
\eeq

Recall the estimate in (\ref{A_pq Decay Est}) and take $\S=\{i\}$. We have (\ref{Characteristic Est2 Case2}) equal to zero at $\lambda=0$. Together with (\ref{Necessity est2}), we find
 \bel{Necessary Est2}
 \begin{array}{rl}\ds
\delta~<~\N_i\left({p-1\over p}\right)+\sum_{j\neq i}\alpha_j-\N_j\left({1\over p}-{1\over q}\right),
\\\\ \ds
i=1,2,\ldots,n.
\end{array}
\eeq
\v

{\bf Case Three:} Consider $\gamma>0, \delta>0$.  Note that (\ref{A-Characteristic sup}) is invariant by changing dilations in one-parameter as shown in  (\ref{Equi Chara}), because of (\ref{Formula}). 

Recall the definition of $\U$ and $\V$ from (\ref{Constraint Case Three}). We write $x_\U\in\R^{\N_\U}$ and $x_\V\in\R^{\N_\V}$ where $\R^{\N_\U}=\bigotimes_{i\in\U}\R^{\N_i}$ and $\R^{\N_\V}=\bigotimes_{i\in\V}\R^{\N_i}$.

Let $|\Q_i|^{1\over \N_i}=\lambda^{-1}$ for every $i\in\U$ and   $|\Q_i|^{1\over \N_i}=1$ for all other  $i\notin\U$. 
We have
\bel{Characteristic Est Case3+}
\begin{array}{lr}\ds
\prod_{i=1}^n |\Q_i|^{{\alpha_i\over \N_i}-\left({1\over p}-{1\over q}\right)}\left\{{1\over|\Q|}\int_\Q \left({1\over|x|}\right)^{\gamma q}dx\right\}^{1\over q}\left\{{1\over|\Q|}\int_\Q \left({1\over |x|}\right)^{\delta p\over p-1}dx\right\}^{p-1\over p}
\\\\ \ds
~\ge~\C_{p~q~\gamma~\delta~n}\prod_{i\in\U}\left({1\over\lambda}\right)^{\alpha_i-{\N_i\over p}}\left\{\idotsint_{\bigotimes_{i\in\U}\Q_i}\left({1\over 1+\sum_{i\in\U}|x_i|}\right)^{\gamma q}\prod_{i\in\U} dx_i\right\}^{1\over q}
\\\\ \ds~~~~~~~~~~~~~~~~~~~~~~~~~~~~~~~~~~~~~~~~~~~
\left\{\prod_{i\in\U}\lambda^{\N_i} \idotsint_{\bigotimes_{i\in\U}\Q_i}\lambda^{\delta p\over p-1}\prod_{ i\in\U}dx_i\right\}^{p-1\over p}\qquad\hbox{\small{($0<\lambda<1$)}}
\\\\ \ds
~\ge~\C_{p~q~\gamma~\delta~n}~\prod_{i\in\U}\left({1\over\lambda}\right)^{\alpha_i-{\N_i\over p}}\left\{\idotsint_{\bigotimes_{i\in\U} 0<|x_i|\leq1}\prod_{i\in\U} dx_i\right\}^{1\over q}
\\\\ \ds~~~~~~~~~~~~~~~~~~~~~~~~~~~~~~~~~~~~~~~~~~~
\left\{\prod_{i\in\U}\lambda^{\N_i} \idotsint_{\bigotimes_{i\in\U}\Q_i}\lambda^{\delta p\over p-1}\prod_{ i\in\U}dx_i\right\}^{p-1\over p}
\\\\ \ds
~\ge~\C_{p~q~\gamma~\delta~n~\N}~\left({1\over \lambda}\right)^{\sum_{i\in\U} \alpha_i-{\N_i\over p}-\delta}.
\end{array}
\eeq
In the case of  $\U=\{1,2,\ldots,n\}$, since $\gamma$ satisfies the first strict inequality in (\ref{local integrability}),  we find 
\bel{Case3 computation1}
\begin{array}{lr}\ds
\delta~=~{\N\over q}-\gamma+\sum_{i=1}^n \alpha_i-{\N_i\over p}\qquad \hbox{by (\ref{Formula})}
\\\\ \ds~~
~>~\sum_{i=1}^n  \alpha_i-{\N_i\over p}~=~\sum_{i\in\U} \alpha_i-{\N_i\over p}.
\end{array}
\eeq
\vsk
Suppose that $\U$ is a proper subset of $\{1,2,\ldots, n\}$ and there exists at least one $i\in\U^c$ such that 
$\alpha_i-\N_i\left({1\over p}-{1\over q}\right)>0$. 
By applying the estimate in (\ref{A_pq Decay Est}) with $\S=\U$, we have (\ref{Characteristic Est Case3+}) equal to zero at $\lambda=0$. 
The last line of (\ref{Characteristic Est Case3+}) implies
\bel{Case3 delta  strict}
\sum_{i\in\U} \alpha_i-{\N_i\over p}~<~\delta.
\eeq

Suppose that $\U$ is a proper subset of $\{1,2,\ldots, n\}$ where 
$\alpha_i-\N_i\left({1\over p}-{1\over q}\right)=0$ for every $i\in\U^c$.

Let $\S=\U$. We have $|x_\U|^{-\gamma}$ and $|x_\U|^\delta$ satisfying the Muckenhoupt characteristic    (\ref{A-subspace S}) on $\R^{\N_\U}\doteq\bigotimes_{i\in\U}\R^{\N_i}$. Denote $\alpha_\U=\sum_{i\in\U}\alpha_i$. By carrying out  the  same estimate in (\ref{Equi Chara}), we find
\bel{Formula U}
\begin{array}{cc}\ds
\gamma~<~{\N_\U\over q},\qquad\delta~<~\N_\U\left({p-1\over p}\right),
\\\\ \ds
 {\alpha_\U\over\N_\U}~=~{1\over p}-{1\over q}+{\gamma+\delta\over\N_\U}. 
 \end{array}
 \eeq
This further implies
\bel{Case3 computation1 U}
\begin{array}{lr}\ds
\delta~=~{\N_\U\over q}-\gamma+\sum_{i\in\U} \alpha_i-{\N_i\over p}
~>~\sum_{i\in\U} \alpha_i-{\N_i\over p}.
\end{array}
\eeq

Let  $|\Q_i|^{1\over \N_i}=\lambda^{-1}$ for every $i\in\V$ and  $|\Q_i|^{1\over \N_i}=1$ for all other  $i\notin\V$. We have
\bel{Characteristic Est Case3-}
\begin{array}{lr}\ds
\prod_{i=1}^n |\Q_i|^{{\alpha_i\over \N_i}-\left({1\over p}-{1\over q}\right)}\left\{{1\over|\Q|}\int_\Q \left({1\over|x|}\right)^{\gamma q}dx\right\}^{1\over q}\left\{{1\over|\Q|}\int_\Q \left({1\over |x|}\right)^{\delta p\over p-1}dx\right\}^{p-1\over p}
\\\\ \ds
~\ge~\C_{p~q~\gamma~\delta~n}~\left({1\over \lambda}\right)^{\sum_{i\in\V}\alpha_i-\left({q-1\over q}\right)\N_i}\left\{\prod_{i\in\V}\lambda^{\N_i}\idotsint_{\bigotimes_{i\in\V}\Q_i}\lambda^{\gamma q}\prod_{i\in\V} dx_i\right\}^{1\over q}
\\\\ \ds~~~~~~~~~~~~~~~~~~~~~~~~~~~~~~~~~~~
~~~~~\left\{ \idotsint_{\bigotimes_{i\in\V}\Q_i}\left({1\over 1+\sum_{ i\in\V}|x_i|}\right)^{\delta p\over p-1}\prod_{ i\in\V}dx_i\right\}^{p-1\over p}
\\\\ \ds
~\ge~\C_{p~q~\gamma~\delta~n}~\left({1\over \lambda}\right)^{\sum_{i\in\V}\alpha_i-\left({q-1\over q}\right)\N_i}\left\{\prod_{i\in\V}\lambda^{\N_i}\idotsint_{\bigotimes_{i\in\V}\Q_i}\lambda^{\gamma q}\prod_{i\in\V} dx_i\right\}^{1\over q}
\\\\ \ds~~~~~~~~~~~~~~~~~~~~~~~~~~~~~~~~~~~

\left\{ \idotsint_{\bigotimes_{i\in\V}0<|x_i|\leq1}\prod_{ i\in\V}dx_i\right\}^{p-1\over p}
\\\\ \ds
~\ge~\C_{p~q~\gamma~\delta~n~\N}~\left({1\over \lambda}\right)^{\sum_{i\in\V} \alpha_i-\left({q-1\over q}\right)\N_i-\gamma}.
\end{array}
\eeq

In the case of $\V=\{1,2,\ldots,n\}$, since $\delta$ satisfies the second strict inequality in (\ref{local integrability}), we find
\bel{Case3 computation2}
\begin{array}{lr}\ds
\gamma~=~\left({p-1\over p}\right)\N-\delta+\sum_{i=1}^n \alpha_i-\N_i\left({q-1\over q}\right)\qquad \hbox{by (\ref{Formula})}
\\\\ \ds~~
~>~\sum_{i=1}^n \alpha_i-\N_i\left({q-1\over q}\right)
~=~ \sum_{i\in\V} \alpha_i-\N_i\left({q-1\over q}\right).
\end{array}
\eeq

Suppose that $\V$ is a proper subset of $\{1,2,\ldots, n\}$ and there exists at least one $i\in\V^c$ such that 
$\alpha_i-\N_i\left({1\over p}-{1\over q}\right)>0$. By applying the estimate in (\ref{A_pq Decay Est}) with $\S=\V$, we have (\ref{Characteristic Est Case3-}) equal to zero at $\lambda=0$. 
The last line of (\ref{Characteristic Est Case3-}) implies
\bel{Case3 gamma strict}
\sum_{i\in\V} \alpha_i-\left({q-1\over q}\right)\N_i~<~\gamma.
\eeq
Suppose that $\V$ is a proper subset of $\{1,2,\ldots, n\}$ where
$\alpha_i-\N_i\left({1\over p}-{1\over q}\right)=0$ for every $i\in\V^c$.

Let $\S=\V$. We have $|x_\V|^{-\gamma}$ and $|x_\V|^\delta$ satisfying the Muckenhoupt characteristic   (\ref{A-subspace S})  on $\R^{\N_\V}\doteq\bigotimes_{i\in\V}\R^{\N_i}$. Denote $\alpha_\V=\sum_{i\in\V}\alpha_i$. By carrying out the same estimate in (\ref{Equi Chara}), we find
\bel{Formula V}
\begin{array}{cc}\ds
\gamma~<~{\N_\V\over q},\qquad\delta~<~\N_\V\left({p-1\over p}\right),
\qquad
 {\alpha_\V\over\N_\V}~=~{1\over p}-{1\over q}+{\gamma+\delta\over\N_\V}.
 \end{array}
 \eeq
This further implies
\bel{Case3 computation2 V}
\begin{array}{lr}\ds
\gamma~=~\left({p-1\over p}\right)\N_\V-\delta+\sum_{i\in\V} \alpha_i-\N_i\left({q-1\over q}\right)
\\\\ \ds~~
~>~\sum_{i\in\V} \alpha_i-\N_i\left({q-1\over q}\right).
\end{array}
\eeq
\begin{remark} By using the formula in (\ref{Formula}), we can verify  that the constraints in (\ref{Necessary Est1}) and (\ref{Necessary Est2}) are equivalent to (\ref{Constraint Case One}) and  (\ref{Constraint Case Two}) respectively. Namely,

for $\gamma\ge0,\delta\leq0$, 
\bel{Necessary Est1 Equiv}
\begin{array}{rl}\ds
\gamma~<~{\N_i\over q}+\sum_{j\neq i} \alpha_j-\N_j\left({1\over p}-{1\over q}\right)\qquad\Longleftrightarrow \qquad\alpha_i-{\N_i\over p}~<~\delta,
\\ \ds
 i~=~1,2,\ldots,n,
 \end{array}
\eeq
 for $\gamma\leq0,\delta\ge0$, 
 \bel{Necessary Est2 Equiv}
 \begin{array}{rl}\ds
\delta~<~\N_i\left({p-1\over p}\right)+\sum_{j\neq i}\alpha_j-\N_j\left({1\over p}-{1\over q}\right)\qquad\Longleftrightarrow\qquad\alpha_i-\N_i\left({q-1\over q}\right)~<~\gamma,~~~~
\\ \ds
 i~=~1,2,\ldots,n.~~~~~~~
 \end{array}
\eeq
\end{remark}

\section{Decay estimate on varying eccentricities}
\setcounter{equation}{0}
{\bf Principal Lemma:~~} {\it Let  $\gamma, \delta$ satisfying  (\ref{local integrability})-(\ref{Constraint Case Three}). Suppose 
\bel{strict subbalance}
{\alpha_i\over \N_i}~>~{1\over p}-{1\over q},\qquad i=1,2,\ldots,n.
\eeq
For $0<\lambda_i\leq1, i=1,2,\ldots,n$,  define
\bel{ratio Q}
{^\lambda}\Q\subset\R^\N~\colon~{|\Q_i|^{1\over \N_i}/|\Q_\nu|^{1\over\N_\nu}}~=~\lambda_i.
\eeq
There exists an $\ve>0$  such that
\bel{Eccentricity Decay}
\begin{array}{rl}\ds
\sup_{\hbox{\small{${^\lambda}\Q$}}}~\prod_{i=1}^n |\Q_i|^{{\alpha_i\over \N_i}-\left({1\over p}-{1\over q}\right)}\left\{{1\over|\Q|}\int_\Q \left({1\over|x|}\right)^{\gamma qr}dx\right\}^{1\over qr}\left\{{1\over|\Q|}\int_\Q \left({1\over |x|}\right)^{\delta pr\over p-1}dx\right\}^{p-1\over pr}
\\\\ \ds
~\leq~ \C_{p~q~r~\alpha~\gamma~\delta~n~\N}~\prod_{i=1}^n \left(\lambda_i\right)^\ve~~~~~~~
\end{array}
\eeq
for some $r>1$. 
The values of $\ve$ and $r$ depend only on $p,q,\gamma, \delta, \alpha, n, \N$. }

\begin{remark} Without the condition (\ref{strict subbalance}), we can only show that the Muckenhoupt characteristic in (\ref{Eccentricity Decay}) is bounded.
\end{remark}

{\bf Proof:}  By carrying out the same estimate in (\ref{Equi Chara}) and using the formula (\ref{Formula}), we find that the {\it $r$-bump} characteristic  (\ref{Eccentricity Decay}) is invariant by changing  dilations in one-parameter. Therefore, it is  suffice to consider $|\Q_\nu|^{1\over\N_\nu}=1$.

Let $\Q_i^o$ and $\Q_i^*\subset\R^{\N_i}$ to be centered on the origin of $\R^{\N_i}$ and 
\bel{Q*}
|\Q_i^o|^{1\over\N_i}~=~|\Q_i|^{1\over\N_i},\qquad |\Q_i^*|^{1\over\N_i}~=~3|\Q_i|^{1\over \N_i}~=~3\lambda_i,\qquad i~=~1,2,\ldots,n.
\eeq
\begin{remark}  
Suppose $\Q_i\cap\Q_i^o=\emptyset$. We must have $|x_i|\ge|x_i^o|/\sqrt{n}$ for every $x_i\in\Q_i$ and every $x_i^o\in\Q_i^o$. Otherwise, if $\Q_i$ intersects $\Q_i^o$, then  $\Q_i\subset\Q_i^*$.
\end{remark}
After a  permutation on indices $i=1,2,\ldots,n$, we can assume $\nu=1$ and
\bel{lambda >}
1~=~\lambda_1~\ge~\lambda_2~\ge~\cdots~\ge~\lambda_n.
\eeq

{\bf Case One:}~~Let $\gamma\ge0, \delta\leq0$ satisfy  (\ref{local integrability})-(\ref{Constraint Case One}).  By adjusting the value of $r$,  we assume
\bel{Range1}
\sum_{i=1}^{m-1}\N_i~<~\gamma qr~<~\sum_{i=1}^m\N_i,\qquad\delta~\leq~0,\qquad 1~\leq m~\leq n.
\eeq

Suppose that $\Q$ is centered on $z\in\R^\N$ for some $|z|\leq3$. We have 
\bel{Decay Est1}
\begin{array}{lr}\ds
\prod_{i=1}^n |\Q_i|^{{\alpha_i\over \N_i}-\left({1\over p}-{1\over q}\right)}\left\{{1\over|\Q|}\int_\Q \left({1\over|x|}\right)^{\gamma qr}dx\right\}^{1\over qr}\left\{{1\over|\Q|}\int_\Q \left({1\over |x|}\right)^{\delta \left({pr\over p-1}\right)}dx\right\}^{p-1\over pr}
\\\\ \ds
~\leq~ \C_{p~q~r~\gamma~\delta~n~\N}\prod_{i=1}^n \left(\lambda_i\right)^{\alpha_i-\N_i\left({1\over p}-{1\over q}\right)}\left\{\prod_{i=1}^n \left({1\over \lambda_i}\right)^{\N_i}\idotsint_{\bigotimes_{i=1}^n \Q_i}  \left({1\over |x_1|+\cdots+|x_n|}\right)^{\gamma qr} dx_1 \cdots dx_n\right\}^{1\over qr}
\\ \ds~~~~~~~~~~~~~~~~~~~~~~~~~~~~~~~~~~~~~~~~~~~~~~~~~~~~~~~~~~~~~~~~~~~~~~~~~~~~~~~~~~~~~~~~~~~~~~~~~~~~~~~~~~~~~~~~~~~~~~~~~~~~~~~~~~~~~~
\hbox{\small{($\delta\leq0$)}}
\\ \ds
~\leq~ \C_{p~q~r~\gamma~\delta~n~\N}\prod_{i=1}^n \left(\lambda_i\right)^{\alpha_i-\N_i\left({1\over p}-{1\over q}\right)}\prod_{i=1}^n\left({1\over \lambda_i}\right)^{\N_i\over qr} 
\\\\ \ds~~~~~~~
\left\{\idotsint_{\bigotimes_{i=m}^n \Q_i}\left\{\idotsint_{\bigotimes_{i=1}^{m-1} \R^{\N_i}}  \left({1\over |x_1|+\cdots+|x_n|}\right)^{\gamma qr} dx_1 \cdots dx_{m-1}\right\}dx_m\cdots dx_n\right\}^{1\over qr}
\\\\ \ds
~\leq~\C_{p~q~r~\gamma~\delta~n~\N} \prod_{i=1}^n \left(\lambda_i\right)^{\alpha_i-\N_i\left({1\over p}-{1\over q}\right)}\prod_{i=1}^n\left({1\over \lambda_i}\right)^{\N_i\over qr} 

\left\{\idotsint_{\bigotimes_{i=m}^n \Q_i} \left({1\over |x_m|+\cdots+|x_{n}|}\right)^{\gamma qr-\sum_{i=1}^{m-1}\N_i } dx_m\cdots dx_n  \right\}^{1\over qr}
\\\\ \ds
~\leq~ \C_{p~q~r~\gamma~\delta~n~\N}\prod_{i=1}^n \left(\lambda_i\right)^{\alpha_i-\N_i\left({1\over p}-{1\over q}\right)}\prod_{i=1}^n\left({1\over \lambda_i}\right)^{\N_i\over qr} 
\left\{\idotsint_{\bigotimes_{i=m}^n \Q_i} \left({1\over |x_m|}\right)^{\gamma qr-\sum_{i=1}^{m-1}\N_i } dx_m\cdots dx_n  \right\}^{1\over qr}
\\\\ \ds
~\leq~\C_{p~q~r~\gamma~\delta~n~\N}\prod_{i=1}^n \left(\lambda_i\right)^{\alpha_i-\N_i\left({1\over p}-{1\over q}\right)}\prod_{i=1}^{m}\left({1\over \lambda_i}\right)^{\N_i\over qr} 
\left\{\int_{ \Q^*_m} \left({1\over |x_m|}\right)^{\gamma qr-\sum_{i=1}^{m-1}\N_i } dx_m  \right\}^{1\over qr}\qquad \hbox{\small{by{\bf Remark 5.2}}}
\\\\ \ds
~\leq~\C_{p~q~r~\gamma~\delta~n~\N}\left(\lambda_m\right)^{{1\over qr}\sum_{i=1}^m \N_i-\gamma}
\prod_{i=1}^{m}\left({1\over \lambda_i}\right)^{\N_i\over qr} \prod_{i=1}^n \left(\lambda_i\right)^{\alpha_i-\N_i\left({1\over p}-{1\over q}\right)}.
\end{array}
\eeq
From direct computation, the formula in the last line of (\ref{Decay Est1}) can be rewritten as
\bel{Case1 Indices Est1}
\begin{array}{lr}\ds
~~~~~~~\left(\lambda_m\right)^{{1\over qr}\sum_{i=1}^m \N_i-\gamma}
\prod_{i=1}^m \left(\lambda_i\right)^{\alpha_i-\N_i\left({1\over p}-{1\over q}\right)-{\N_i\over qr}}\prod_{i=m+1}^n \left(\lambda_i\right)^{\alpha_i-\N_i\left({1\over p}-{1\over q}\right)}
\\\\ \ds 
~=~\left(\lambda_m\right)^{{1\over qr}\sum_{i=1}^m \N_i-\gamma}
\prod_{i=2}^m \left(\lambda_i\right)^{\alpha_i-\N_i\left({1\over p}-{1\over q}\right)-{\N_i\over qr}}\prod_{i=m+1}^n \left(\lambda_i\right)^{\alpha_i-\N_i\left({1\over p}-{1\over q}\right)}\qquad (\lambda_1=1)
\\\\ \ds
~=~\left(\lambda_m\right)^{{\N_1\over qr}+\sum_{i=2}^n \alpha_i-\N_i\left({1\over p}-{1\over q}\right)-\gamma}
\prod_{i=2}^m \left({\lambda_i\over \lambda_m}\right)^{\alpha_i-{\N_i\over p}+\left(1-{1\over r}\right){\N_i\over q}}\prod_{i=m+1}^n \left({\lambda_i\over \lambda_m}\right)^{\alpha_i-\N_i\left({1\over p}-{1\over q}\right)}.
\end{array}
\eeq
Recall {\bf Remark 4.1}.  $\gamma\ge0,\delta\leq0$ satisfy the two equivalent strict inequalities in (\ref{Necessary Est1 Equiv}).

Define $0\leq\vartheta\leq1$ implicitly by  letting $\lambda_m=(\lambda_n)^\vartheta$.   For $r$ sufficiently close to $1$, we have
\bel{Case1 Indices Est2}
\begin{array}{lr}\ds
\vartheta\left[{\N_1\over qr}+\sum_{i=2}^{n} \alpha_i-\N_i\left({1\over p}-{1\over q}\right)-\gamma\right]+(1-\vartheta)\left[\alpha_n-\N_n\left({1\over p}-{1\over q}\right)\right]~>~0
\end{array}
\eeq
and
\bel{Case1 Indices Est3}
\alpha_i-{\N_i\over p}+\left(1-{1\over r}\right){\N_i \over q}~<~0,\qquad i=1,2,\ldots n.
\eeq
Note that $\alpha_n-\N_n\left({1\over p}-{1\over q}\right)>0$.
By using  (\ref{Case1 Indices Est2})-(\ref{Case1 Indices Est3}), we find that (\ref{Case1 Indices Est1}), is bounded by $\C_{p~q~r~\gamma~\delta~n~\N} ~(\lambda_n)^\ve$ for some $\ve=\ve(p,q,r,\alpha,\gamma,\delta,n,\N)>0$.

{\bf Case Two:~} Let $\gamma\leq0, \delta\ge0$ satisfy  (\ref{local integrability})-(\ref{Formula}) and (\ref{Constraint Case Two}).  By adjusting the value of $r$,  assume
\bel{Range2}
\gamma~\leq~0,\qquad\sum_{i=1}^{m-1}\N_i~<~\delta \left({pr\over p-1}\right)~<~\sum_{i=1}^m\N_i,\qquad 1~\leq~m~\leq~n.
\eeq 
Suppose that $\Q$ is centered on $z\in\R^\N$ for some $|z|\leq3$.
We have 
\bel{Decay Est2}
\begin{array}{lr}\ds
\prod_{i=1}^n |\Q_i|^{{\alpha_i\over \N_i}-\left({1\over p}-{1\over q}\right)}\left\{{1\over|\Q|}\int_\Q \left({1\over|x|}\right)^{\gamma qr}dx\right\}^{1\over qr}\left\{{1\over|\Q|}\int_\Q \left({1\over |x|}\right)^{\delta \left({pr\over p-1}\right)}dx\right\}^{p-1\over pr}
\\\\ \ds
~\leq~\C_{p~q~r~\gamma~\delta~n~\N} \prod_{i=1}^n \left(\lambda_i\right)^{\alpha_i-\N_i\left({1\over p}-{1\over q}\right)}\left\{\prod_{i=1}^n\left({1\over \lambda_i}\right)^{\N_i}\idotsint_{\bigotimes_{i=1}^n \Q_i}  \left({1\over |x_1|+\cdots+|x_n|}\right)^{\delta\left({pr\over p-1}\right)} dx_1 \cdots dx_n\right\}^{p-1\over pr}
\\ \ds~~~~~~~~~~~~~~~~~~~~~~~~~~~~~~~~~~~~~~~~~~~~~~~~~~~~~~~~~~~~~~~~~~~~~~~~~~~~~~~~~~~~~~~~~~~~~~~~~~~~~~~~~~~~~~~~~~~~~~~~~~~~~~~~~~~~~~~~~~~~~~~
\hbox{\small{($\gamma\leq0$)}}
\\ \ds
~\leq~\C_{p~q~r~\gamma~\delta~n~\N}\prod_{i=1}^n \left(\lambda_i\right)^{\alpha_i-\N_i\left({1\over p}-{1\over q}\right)}\prod_{i=1}^n \left({1\over \lambda_i}\right)^{\left({p-1\over pr}\right)\N_i}
\\\\ \ds~~~~~~~
\left\{\idotsint_{\bigotimes_{i=m}^n \Q_i}\left\{\idotsint_{\bigotimes_{i=1}^{m-1}\R^{\N_i}}  \left({1\over |x_1|+\cdots+|x_n|}\right)^{\delta\left({pr\over p-1}\right)} dx_1 \cdots dx_{m-1}\right\}dx_m\cdots dx_n\right\}^{p-1\over pr}
\\\\ \ds
~\leq~ \C_{p~q~r~\gamma~\delta~n~\N}\prod_{i=1}^n \left(\lambda_i\right)^{\alpha_i-\N_i\left({1\over p}-{1\over q}\right)}\prod_{i=1}^n \left({1\over \lambda_i}\right)^{\left({p-1\over pr}\right)\N_i}
\\\\ \ds~~~~~~~
\left\{\idotsint_{\bigotimes_{i=m}^n \Q_i} \left({1\over |x_m|+\cdots+|x_n|}\right)^{\delta\left({pr\over p-1}\right)-\sum_{i=1}^{m-1}\N_i} dx_m\cdots dx_n\right\}^{p-1\over pr}
\\\\ \ds
~\leq~\C_{p~q~r~\gamma~\delta~n~\N} \prod_{i=1}^n \left(\lambda_i\right)^{\alpha_i-\N_i\left({1\over p}-{1\over q}\right)}\prod_{i=1}^n \left({1\over \lambda_i}\right)^{\left({p-1\over pr}\right)\N_i}

\left\{\idotsint_{\bigotimes_{i=m}^n \Q_i} \left({1\over |x_m|}\right)^{\delta\left({pr\over p-1}\right)-\sum_{i=1}^{m-1}\N_i} dx_m\cdots dx_n\right\}^{p-1\over pr}

\\\\ \ds
~\leq~ \C_{p~q~r~\gamma~\delta~n~\N}\prod_{i=1}^n \left(\lambda_i\right)^{\alpha_i-\N_i\left({1\over p}-{1\over q}\right)}
\prod_{i=1}^m \left({1\over \lambda_i}\right)^{\left({p-1\over pr}\right)\N_i}
\left\{ \int_{\Q^*_m} \left({1\over |x_{m}|}\right)^{\delta\left({pr\over p-1}\right)-\sum_{i=1}^{m-1}\N_i }dx_m \right\}^{p-1\over pr}~~\hbox{\small{by {\bf Remark 5.2}}}
\\\\ \ds
~\leq~\C_{p~q~r~\gamma~\delta~n~\N}
\left(\lambda_m\right)^{\left({p-1\over pr}\right)\sum_{i=1}^m\N_i -\delta} 
\prod_{i=1}^m \left({1\over \lambda_i}\right)^{\left({p-1\over pr}\right)\N_i}\prod_{i=1}^n \left(\lambda_i\right)^{\alpha_i-\N_i\left({1\over p}-{1\over q}\right)}.
\end{array}
\eeq

From direct computation, the formula in the last line of (\ref{Decay Est2}) can be rewritten as
\bel{Case2 Indices Est1}
\begin{array}{lr}\ds
~~~~~~~\left(\lambda_m\right)^{\left({p-1\over pr}\right)\sum_{i=1}^m\N_i -\delta} 
\prod_{i=1}^m \left(\lambda_i\right)^{\alpha_i-\N_i\left({1\over p}-{1\over q}\right)-\left({p-1\over pr}\right)\N_i}
\prod_{i=m+1}^n \left( \lambda_i\right)^{\alpha_i-\N_i\left({1\over p}-{1\over q}\right)}
\\\\ \ds
~=~\left(\lambda_m\right)^{\left({p-1\over pr}\right)\sum_{i=1}^m \N_i-\delta} 

\prod_{i=2}^m \left(\lambda_i\right)^{\alpha_i-\N_i\left({1\over p}-{1\over q}\right)-\left({p-1\over pr}\right)\N_i}
\prod_{i=m+1}^n \left( \lambda_i\right)^{\alpha_i-\N_i\left({1\over p}-{1\over q}\right)}\qquad (\lambda_1=1)
\\\\ \ds
~=~\left(\lambda_m\right)^{\left({p-1\over pr}\right)\N_1+\sum_{i=2}^n \alpha_i-\N_i\left({1\over p}-{1\over q}\right)-\delta} 

\prod_{i=2}^m \left({\lambda_i\over \lambda_m}\right)^{\alpha_i-\N_i\left({q-1\over q}\right)+\left(1-{1\over r}\right)\left({p-1\over p}\right)\N_i}
\prod_{i=m+1}^n \left( {\lambda_i\over \lambda_m}\right)^{\alpha_i-\N_i\left({1\over p}-{1\over q}\right)}.
\end{array}
\eeq
Recall {\bf Remark 4.1}.  $\gamma\leq0,\delta\ge0$ satisfy the two equivalent strict inequalities in (\ref{Necessary Est2 Equiv}).  
Define $0\leq\vartheta\leq1$ implicitly by letting $\lambda_m=(\lambda_n)^\vartheta$. 
For $r$ sufficiently close to $1$, we have 
 \bel{Case2 Indices Est2}
\begin{array}{lr}\ds
\vartheta\left[\N_1\left({p-1\over pr}\right)+\sum_{i=2}^{n} \alpha_i-\N_i\left({1\over p}-{1\over q}\right)-\delta\right]+(1-\vartheta)\left[\alpha_n-\N_n\left({1\over p}-{1\over q}\right)\right]~>~0
\end{array}
\eeq
and
\bel{Case2 Indices Est3}
\alpha_i-\N_i\left({q-1\over q}\right)+\left(1-{1\over r}\right)\left({p-1\over p}\right)\N_i~ <~0,\qquad i=1,2,\ldots n.
\eeq
Note that $\alpha_n-\N_n\left({1\over p}-{1\over q}\right)>0$. By using (\ref{Case2 Indices Est2})-(\ref{Case2 Indices Est3}), we find that  (\ref{Case2 Indices Est1}) is  bounded by  $\C_{p~q~r~\gamma~\delta~n~\N} ~(\lambda_n)^\ve$ for some $\ve=\ve(p,q,r,\alpha,\gamma,\delta,n,\N)>0$.

Suppose that  $\Q$ is centered on $z\in\R^\N$ for which $|z|>3$. Since $\Q$ has a diameter $1$,  we have 
\bel{distance inclusion}
{1\over 2}|z|~\leq~|x|~\leq~2|z|
\eeq
whenever $x\in\Q$. From (\ref{distance inclusion}), we  have
\bel{Decay est}
\begin{array}{lr}\ds
~~~~~~~\prod_{i=1}^n |\Q_i|^{{\alpha_i\over \N_i}-\left({1\over p}-{1\over q}\right)}\left\{{1\over|\Q|}\int_\Q \left({1\over|x|}\right)^{\gamma qr}dx\right\}^{1\over qr}\left\{{1\over|\Q|}\int_\Q \left({1\over |x|}\right)^{\delta\left({pr\over p-1}\right)}dx\right\}^{p-1\over pr}
\\\\ \ds
~\leq~\C_{\gamma~\delta}\left({1\over|z|}\right)^{\gamma+\delta}\prod_{i=1}^n \left(\lambda_i\right)^{\alpha_i-\N_i\left({1\over p}-{1\over q}\right)}~\leq~\C_{\gamma~\delta}\prod_{i=1}^n \left(\lambda_i\right)^{\alpha_i-\N_i\left({1\over p}-{1\over q}\right)}.\qquad\hbox{\small{($\gamma+\delta\ge0$)}}
\end{array}
 \eeq
{\bf Case Three:~~} Let $\gamma>0, \delta>0$ satisfy  (\ref{local integrability})-(\ref{Formula}) and (\ref{Constraint Case Three}).  By adjusting the value of $r$,  assume
\bel{Range3 gamma}
 \sum_{i=1}^{m-1}\N_i~<~\gamma qr~<~\sum_{i=1}^m\N_i,\qquad 1~\leq~m~\leq~n,
 \eeq
 \bel{Range3 delta}
  \sum_{i=1}^{l-1}\N_i~<~\delta \left({pr\over p-1}\right)~<~\sum_{i=1}^l\N_i,\qquad 1~\leq ~l~\leq~n.
  \eeq
We have
\bel{Decay Case3 Est}
\begin{array}{lr}\ds
\prod_{i=1}^n |\Q_i|^{{\alpha_i\over \N_i}-\left({1\over p}-{1\over q}\right)}\left\{{1\over|\Q|}\int_\Q \left({1\over|x|}\right)^{\gamma qr}dx\right\}^{1\over qr}\left\{{1\over|\Q|}\int_\Q \left({1\over |x|}\right)^{\delta\left({pr\over p-1}\right)}dx\right\}^{p-1\over pr}
\\\\ \ds
~\leq~\C_{p~q~r~\gamma~\delta~n} \prod_{i=1}^n\left(\lambda_i\right)^{\alpha_i-\N_i\left({1\over p}-{1\over q}\right)}\prod_{i=1}^n\left({1\over \lambda_i}\right)^{\N_i\over qr} \prod_{i=1}^n \left({1\over \lambda_i}\right)^{\N_i\left({p-1\over pr}\right)}
\\\\ \ds~~~~~~~
\left\{\idotsint_{\bigotimes_{i=m}^n\Q_i}\left\{\idotsint_{\bigotimes_{i=1}^{m-1} \R^{\N_i}}  \left({1\over |x_1|+\cdots+|x_n|}\right)^{\gamma qr} dx_1 \cdots dx_{m-1}\right\} dx_m\cdots dx_n\right\}^{1\over qr}
\\\\ \ds~~~~~~~
\left\{\idotsint_{\bigotimes_{i=l}^n \Q_i}\left\{\idotsint_{\bigotimes_{i=1}^{l-1} \R^{\N_i}}  \left({1\over |x_1|+\cdots+|x_n|}\right)^{\delta\left({pr\over p-1}\right)} dx_1 \cdots dx_{l-1}\right\}dx_l\cdots dx_n\right\}^{p-1\over pr}
\\\\ \ds
~\leq~\C_{p~q~r~\gamma~\delta~n~\N} \prod_{i=1}^n\left(\lambda_i\right)^{\alpha_i-\N_i\left({1\over p}-{1\over q}\right)}\prod_{i=1}^n\left({1\over \lambda_i}\right)^{\N_i\over qr} \prod_{i=1}^n \left({1\over \lambda_i}\right)^{\N_i\left({p-1\over pr}\right)}
\\\\ \ds~~~~~~~
\left\{\idotsint_{\bigotimes_{i=m}^n\Q_i} \left({1\over |x_m|+\cdots+|x_n|}\right)^{\gamma qr-\sum_{i=1}^{m-1}\N_i}  dx_m\cdots dx_n\right\}^{1\over qr}
\\\\ \ds~~~~~~~
\left\{\idotsint_{\bigotimes_{i=l}^n \Q_i}  \left({1\over |x_l|+\cdots+|x_n|}\right)^{\delta\left({pr\over p-1}\right)-\sum_{i=1}^{l-1}\N_i} dx_l\cdots dx_n\right\}^{p-1\over pr}
\\\\ \ds
~\leq~\C_{p~q~r~\gamma~\delta~n~\N} \prod_{i=1}^n\left(\lambda_i\right)^{\alpha_i-\N_i\left({1\over p}-{1\over q}\right)}\prod_{i=1}^n\left({1\over \lambda_i}\right)^{\N_i\over qr} \prod_{i=1}^n \left({1\over \lambda_i}\right)^{\N_i\left({p-1\over pr}\right)}
\\\\ \ds~~~~~~~
\left\{\idotsint_{\bigotimes_{i=m}^n\Q_i} \left({1\over |x_m|}\right)^{\gamma qr-\sum_{i=1}^{m-1}\N_i}  dx_m\cdots dx_n\right\}^{1\over qr}
\left\{\idotsint_{\bigotimes_{i=l}^n \Q_i}  \left({1\over |x_l|}\right)^{\delta\left({pr\over p-1}\right)-\sum_{i=1}^{l-1}\N_i} dx_l\cdots dx_n\right\}^{p-1\over pr}
\\\\ \ds
~\leq~\C_{p~q~r~\gamma~\delta~n~\N}\prod_{i=1}^n\left(\lambda_i\right)^{\alpha_i-\N_i\left({1\over p}-{1\over q}\right)}\prod_{i=1}^m\left({1\over \lambda_i}\right)^{\N_i\over qr} \prod_{i=1}^l \left({1\over \lambda_i}\right)^{\N_i\left({p-1\over pr}\right)}
\\\\ \ds~~~~~~~
\left\{\int_{ \Q^*_m}  \left({1\over |x_m|}\right)^{\gamma qr-\sum_{i=1}^{m-1}\N_i}   dx_m\right\}^{1\over qr}
\left\{\int_{ \Q^*_l}  \left({1\over |x_l|}\right)
^{\delta\left( {pr\over p-1}\right)-\sum_{i=1}^{l-1}\N_i}   dx_l\right\}^{p-1\over pr}
\qquad \hbox{\small{by {\bf Remark 5.2}}}
\\\\ \ds
~\leq~\C_{p~q~r~\gamma~\delta~n~\N}~ \left(\lambda_m\right)^{{1\over qr}\sum_{i=1}^m \N_i-\gamma}\left(\lambda_l\right)^{\left({p-1\over pr}\right)\sum_{i=1}^l\N_i-\delta}\prod_{i=1}^n\left(\lambda_i\right)^{\alpha_i-\N_i\left({1\over p}-{1\over q}\right)}\prod_{i=1}^m\left({1\over \lambda_i}\right)^{\N_i\over qr} \prod_{i=1}^l \left({1\over \lambda_i}\right)^{\left({p-1\over pr}\right)\N_i}.
\end{array}
\eeq

Let $0\leq k\leq n-1$. From direct computation,  we have
\bel{crucial est}
\begin{array}{lr}\ds
{1\over r}\left({1\over q}+{p-1\over p}\right) \sum_{i=1}^k \N_i -(\gamma+\delta)+ \sum_{i=k+1}^n \alpha_i-\N_i\left({1\over p}-{1\over q}\right)
\\\\ \ds
~=~{\N\over r}-{1\over r}\left({1\over p}-{1\over q}\right)\N-(\gamma+\delta)+ \sum_{i=k+1}^n \alpha_i-{\N_i\over r} -\N_i\left(1-{1\over r}\right)\left({1\over p}-{1\over q}\right)
\\\\ \ds
~=~{\N\over r}-{1\over r}\left({1\over p}-{1\over q}\right)\N-\alpha+\N\left({1\over p}-{1\over q}\right)+ \sum_{i=k+1}^n \alpha_i-{\N_i\over r} -\N_i\left(1-{1\over r}\right)\left({1\over p}-{1\over q}\right)\qquad \hbox{by (\ref{Formula})}
\\\\ \ds
~=~\left({\N\over r}-\alpha\right)+\N\left(1-{1\over r}\right)\left({1\over p}-{1\over q}\right)+ \sum_{i=k+1}^n \alpha_i-{\N_i\over r} -\N_i\left(1-{1\over r}\right)\left({1\over p}-{1\over q}\right)
\\\\ \ds
~=~\sum_{i=1}^k {\N_i\over r}-\alpha_i+\N_i\left(1-{1\over r}\right)\left({1\over p}-{1\over q}\right).
\end{array}
\eeq
Suppose $ l\leq m$. The formula in the last line of (\ref{Decay Case3 Est}) can be rewritten as 
\bel{Case3 Indices Est1+}
\begin{array}{lr}\ds
\left(\lambda_m\right)^{{1\over r}\left({1\over q}+{p-1\over p}\right)\sum_{i=1}^{l-1} \N_i-(\gamma+\delta)}\left({\lambda_l\over \lambda_m}\right)^{\left({p-1\over pr}\right)\sum_{i=1}^{l-1}\N_i-\delta}

\prod_{i=1}^n\left(\lambda_i\right)^{\alpha_i-\N_i\left({1\over p}-{1\over q}\right)}\prod_{i=l}^m\left({\lambda_m\over \lambda_i}\right)^{\N_i\over qr} \prod_{i=1}^{l-1} \left({1\over \lambda_i}\right)^{{1\over r}\left({1\over q}+{p-1\over p}\right)\N_i}
\\\\ \ds
~=~\left(\lambda_m\right)^{{1\over r}\left({1\over q}+{p-1\over p}\right)\sum_{i=1}^{l-1} \N_i+\sum_{i=l}^n \alpha_i-\N_i\left({1\over p}-{1\over q}\right)  -(\gamma+\delta)}\left({\lambda_l\over \lambda_m}\right)^{\left({p-1\over pr}\right)\sum_{i=1}^{l-1}\N_i-\delta}
\\\\ \ds~~~~~~~

\prod_{i=l}^n\left({\lambda_i\over \lambda_m}\right)^{\alpha_i-\N_i\left({1\over p}-{1\over q}\right)}
\prod_{i=l}^m\left({\lambda_m\over \lambda_i}\right)^{\N_i\over qr} 
\prod_{i=1}^{l-1} \left({1\over \lambda_i}\right)^{{\N_i\over r}-\alpha_i+\left(1-{1\over r}\right)\left({1\over p}-{1\over q}\right)\N_i}
\\\\ \ds
~=~\left({\lambda_l\over \lambda_m}\right)^{\left({p-1\over pr}\right)\sum_{i=1}^{l-1}\N_i-\delta}
\prod_{i=l}^n\left({\lambda_i\over \lambda_m}\right)^{\alpha_i-\N_i\left({1\over p}-{1\over q}\right)}
\prod_{i=l}^m\left({\lambda_m\over \lambda_i}\right)^{\N_i\over qr} 
\prod_{i=1}^{l-1} \left({\lambda_m\over \lambda_i}\right)^{{\N_i\over r}-\alpha_i+\left(1-{1\over r}\right)\left({1\over p}-{1\over q}\right)\N_i}
\qquad \hbox{by (\ref{crucial est})}
\\\\ \ds
~=~\left({\lambda_l\over \lambda_m}\right)^{\left({p-1\over pr}\right)\sum_{i=1}^{l-1}\N_i-\delta}
 \prod_{i=m+1}^n\left({\lambda_i\over \lambda_m}\right)^{\alpha_i-\N_i\left({1\over p}-{1\over q}\right)}
\prod_{i=l}^{m} \left({\lambda_i\over \lambda_m}\right)^{\alpha_i-{\N_i\over p}+\left(1-{1\over r}\right){\N_i\over q}}
\prod_{i=1}^{l-1} \left({\lambda_m\over \lambda_i}\right)^{{\N_i\over r}-\alpha_i+\left(1-{1\over r}\right)\left({1\over p}-{1\over q}\right)\N_i}
\\\\ \ds
~=~\prod_{i=m+1}^n\left({\lambda_i\over \lambda_m}\right)^{\alpha_i-\N_i\left({1\over p}-{1\over q}\right)}
\prod_{i=1}^{l-1} \left({\lambda_l\over \lambda_i}\right)^{{\N_i\over r}-\alpha_i+\left(1-{1\over r}\right)\left({1\over p}-{1\over q}\right)\N_i}
\\\\ \ds~~~~~~~
\left({\lambda_m\over \lambda_l}\right)^{\sum_{i=1}^{l-1}{\N_i\over p}-\alpha_i-\left(1-{1\over r}\right){\N_i\over q}+\delta}\prod_{i=l}^{m} \left({\lambda_m\over \lambda_i}\right)^{{\N_i\over p}-\alpha_i-\left(1-{1\over r}\right){\N_i\over q}}.
\end{array}
\eeq
Recall the subset $\U$ defined in (\ref{Constraint Case Three}) where $\alpha_i-\N_i/p<0$ for every $i\notin\U$. 

Notice that $\lambda_m\leq\lambda_l$ when $l\leq m$.
For $r$ sufficiently close to $1$, 
we have
 \bel{Case3 Indices Est2+}
 \begin{array}{lr}\ds
 \left({\lambda_m\over \lambda_l}\right)^{\sum_{i=1}^{l-1}{\N_i\over p}-\alpha_i-\left(1-{1\over r}\right){\N_i\over q}+\delta}
\prod_{i=l}^{m} \left({\lambda_m\over \lambda_i}\right)^{{\N_i\over p}-\alpha_i-\left(1-{1\over r}\right){\N_i\over q}}
\\\\ \ds
~\leq~ \left({\lambda_m\over \lambda_l}\right)^{\sum_{i\in\U\cap\{1,\ldots,l-1\}}{\N_i\over p}-\alpha_i-\left(1-{1\over r}\right){\N_i\over q}+\delta}
  \left({\lambda_m\over\lambda_l}\right)^{\sum_{i\in\U\cap\{l,\ldots,m\}} {\N_i\over p}-\alpha_i-\left(1-{1\over r}\right){\N_i\over q}}
 \\\\ \ds~~~~~~~

  \left({\lambda_m\over \lambda_l}\right)^{\sum_{i\in\U^c\cup\{1,\ldots,l-1\}}{\N_i\over p}-\alpha_i-\left(1-{1\over r}\right){\N_i\over q}} 
  \prod_{i\in\U^c\cap\{l,\ldots,m\}} \left({\lambda_m\over \lambda_i}\right)^{{\N_i\over p}-\alpha_i-\left(1-{1\over r}\right){\N_i\over q}}
  \\\\ \ds
 ~\leq~ \left({\lambda_m\over\lambda_l}\right)^{\sum_{i\in\U\cup\{1,\ldots,m\}} {\N_i\over p}-\alpha_i-\left(1-{1\over r}\right){\N_i\over q}+\delta}.
  \end{array}
 \eeq
By bringing the estimates in (\ref{Case3 Indices Est1+})-(\ref{Case3 Indices Est2+}) back to (\ref{Decay Case3  Est}), we find
 \bel{Case3 Indices Result+}
 \begin{array}{lr}\ds
 \left(\lambda_m\right)^{{1\over qr}\sum_{i=1}^m \N_i-\gamma}\left(\lambda_l\right)^{\left({p-1\over pr}\right)\sum_{i=1}^l\N_i-\delta}\prod_{i=1}^n\left(\lambda_i\right)^{\alpha_i-\N_i\left({1\over p}-{1\over q}\right)} 
 \prod_{i=1}^m\left({1\over \lambda_i}\right)^{\N_i\over qr} \prod_{i=1}^l \left({1\over \lambda_i}\right)^{\left({p-1\over pr}\right)\N_i}
 \\\\ \ds
 ~=~\left({\lambda_m\over \lambda_l}\right)^{\sum_{i=1}^{l-1}{\N_i\over p}-\alpha_i-\left(1-{1\over r}\right){\N_i\over q}+\delta}\prod_{i=l}^{m} \left({\lambda_m\over \lambda_i}\right)^{{\N_i\over p}-\alpha_i-\left(1-{1\over r}\right){\N_i\over q}} 
 \\\\ \ds
~~~~~~~\prod_{i=m+1}^n\left({\lambda_i\over \lambda_m}\right)^{\alpha_i-\N_i\left({1\over p}-{1\over q}\right)}
\prod_{i=1}^{l-1} \left({\lambda_l\over \lambda_i}\right)^{{\N_i\over r}-\alpha_i+\left(1-{1\over r}\right)\left({1\over p}-{1\over q}\right)\N_i} 
 \\\\ \ds
~\leq~ 
\left({\lambda_m\over\lambda_l}\right)^{\sum_{i\in\U\cup\{1,2,\ldots,m\}} {\N_i\over p}-\alpha_i-\left(1-{1\over r}\right){\N_i\over q}+\delta}
\prod_{i=1}^{l} \left({\lambda_l\over \lambda_i}\right)^{{\N_i\over r}-\alpha_i+\left(1-{1\over r}\right)\left({1\over p}-{1\over q}\right)\N_i}
\prod_{i=m+1}^n\left({\lambda_i\over \lambda_m}\right)^{\alpha_i-\N_i\left({1\over p}-{1\over q}\right)}.
\end{array}
\eeq
Recall that $\delta>0$ satisfies the first strict inequality in (\ref{Constraint Case Three}).  From (\ref{strict subbalance}) and (\ref{alpha}), we also have $\left({1\over p}-{1\over q}\right)\N_i<\alpha_i<\N_i $ for every $i=1,2,\ldots,n$. Define implicitly $0\leq\vartheta_1\leq\vartheta_2\leq1$ by letting $\lambda_l=(\lambda_n)^{\vartheta_1}$ and $\lambda_m=(\lambda_n)^{\vartheta_2}$.

For $r$ sufficiently close to $1$, we have
\bel{Case3+ Indices Est}
\begin{array}{lr}\ds
~~~~~~~\vartheta_1\left[ {\N_1\over r}-\alpha_1+\left(1-{1\over r}\right)\left({1\over p}-{1\over q}\right)\N_1 \right]
~+~(1-\vartheta_2)\left( \alpha_n-\N_n\left({1\over p}-{1\over q}  \right)   \right)  
\\\\ \ds~~~~~~~
~+~
(\vartheta_2-\vartheta_1)\left[\sum_{i\in\U\cup\{1,2,\ldots,m\}} {\N_i\over p}-\alpha_i-\left(1-{1\over r}\right){\N_i\over q}+\delta\right]
~>~0.
\end{array}
\eeq
The estimate in (\ref{Case3+ Indices Est}) implies that
(\ref{Case3 Indices Result+}) is bounded by a constant multiple of $(\lambda_n)^\ve$ for some $\ve=\ve(p~q~r~\alpha~\gamma~\delta~n~\N)>0$.

On the other hand, suppose $ m\leq l$. The last line of (\ref{Decay Case3 Est}) can be rewritten as 
\bel{Case3 Indices Est1-}
\begin{array}{lr}\ds
\left(\lambda_m\right)^{{1\over qr}\sum_{i=1}^{m-1} \N_i-\gamma}\left(\lambda_l\right)^{\left({p-1\over pr}\right)\sum_{i=1}^{l}\N_i-\delta}\prod_{i=1}^n\left(\lambda_i\right)^{\alpha_i-\N_i\left({1\over p}-{1\over q}\right)}\prod_{i=1}^{m-1}\left({1\over \lambda_i}\right)^{\N_i\over qr} \prod_{i=1}^{l} \left({1\over \lambda_i}\right)^{\left({p-1\over pr}\right)\N_i}
\\\\ \ds
~=~\left(\lambda_l\right)^{{1\over r}\left({1\over q}+{p-1\over p}\right)\sum_{i=1}^{m-1} \N_i-(\gamma+\delta)}\left({\lambda_m\over \lambda_l}\right)^{{1\over qr}\sum_{i=1}^{m-1}\N_i-\gamma}
\prod_{i=1}^n\left(\lambda_i\right)^{\alpha_i-\N_i\left({1\over p}-{1\over q}\right)}\prod_{i=m}^l\left({\lambda_l\over \lambda_i}\right)^{\left({p-1\over pr}\right)\N_i} \prod_{i=1}^{m-1} \left({1\over \lambda_i}\right)^{{1\over r}\left({1\over q}+{p-1\over p}\right)\N_i}
\\\\ \ds
~=~\left(\lambda_l\right)^{{1\over r}\left({1\over q}+{p-1\over p}\right)\sum_{i=1}^{m-1} \N_i+\sum_{i=m}^n \alpha_i-\N_i\left({1\over p}-{1\over q}\right)  -(\gamma+\delta)}\left({\lambda_m\over \lambda_l}\right)^{{1\over qr}\sum_{i=1}^{m-1}\N_i-\gamma}
\\\\ \ds~~~~~~~

\prod_{i=m}^n\left({\lambda_i\over \lambda_l}\right)^{\alpha_i-\N_i\left({1\over p}-{1\over q}\right)}
\prod_{i=m}^l\left({\lambda_l\over \lambda_i}\right)^{\left({p-1\over pr}\right)\N_i} 
\prod_{i=1}^{m-1} \left({1\over \lambda_i}\right)^{{\N_i\over r}-\alpha_i+\left(1-{1\over r}\right)\left({1\over p}-{1\over q}\right)\N_i}
\\\\ \ds
~=~\left({\lambda_m\over \lambda_l}\right)^{{1\over qr}\sum_{i=1}^{m-1}\N_i-\gamma}
\prod_{i=m}^n\left({\lambda_i\over \lambda_l}\right)^{\alpha_i-\N_i\left({1\over p}-{1\over q}\right)}
\prod_{i=m}^l\left({\lambda_l\over \lambda_i}\right)^{\left({p-1\over pr}\right)\N_i}
\prod_{i=1}^{m-1} \left({\lambda_l\over \lambda_i}\right)^{{\N_i\over r}-\alpha_i+\left(1-{1\over r}\right)\left({1\over p}-{1\over q}\right)\N_i} 
\qquad \hbox{by (\ref{crucial est})}
\\\\ \ds
~=~\left({\lambda_m\over \lambda_l}\right)^{{1\over qr}\sum_{i=1}^{m-1}\N_i-\gamma}\prod_{i=l+1}^n\left({\lambda_i\over \lambda_l}\right)^{\alpha_i-\N_i\left({1\over p}-{1\over q}\right)}
\prod_{i=m}^{l} \left({\lambda_i\over \lambda_l}\right)^{\alpha_i-\left({q-1\over q}\right)\N_i+\left(1-{1\over r}\right)\left({p-1\over p}\right)\N_i} 
\prod_{i=1}^{m-1} \left({\lambda_l\over \lambda_i}\right)^{{\N_i\over r}-\alpha_i+\left(1-{1\over r}\right)\left({1\over p}-{1\over q}\right)\N_i}
\\\\ \ds
~=~\prod_{i=l+1}^n\left({\lambda_i\over \lambda_l}\right)^{\alpha_i-\N_i\left({1\over p}-{1\over q}\right)}\prod_{i=1}^{m} \left({\lambda_m\over \lambda_i}\right)^{{\N_i\over r}-\alpha_i+\left(1-{1\over r}\right)\left({1\over p}-{1\over q}\right)\N_i}
\\\\ \ds~~~~~~~
\left({\lambda_l\over \lambda_m}\right)^{\sum_{i=1}^{m-1}\left({q-1\over q}\right)\N_i-\alpha_i-\left(1-{1\over r}\right)\left({p-1\over p}\right)\N_i+\gamma}
\prod_{i=m}^{l} \left({\lambda_l\over \lambda_i}\right)^{\left({q-1\over q}\right)\N_i-\alpha_i-\left(1-{1\over r}\right)\left({p-1\over p}\right)\N_i}.
\end{array}
\eeq
Recall the subset $\V$  defined in (\ref{Constraint Case Three}) where $\alpha_i-\N_i\left({q-1\over q}\right)<0$ for every $i\notin\V$. 

Notice that $\lambda_l\leq\lambda_m$ when $m\leq l$.
For $r$ sufficiently close to $1$, 
we have
 \bel{Case3 Indices Est2-}
 \begin{array}{lr}\ds
\left({\lambda_l\over \lambda_m}\right)^{\sum_{i=1}^{m-1}\left({q-1\over q}\right)\N_i-\alpha_i-\left(1-{1\over r}\right)\left({p-1\over p}\right)\N_i+\gamma}
\prod_{i=m}^{l} \left({\lambda_l\over \lambda_i}\right)^{\left({q-1\over q}\right)\N_i-\alpha_i-\left(1-{1\over r}\right)\left({p-1\over p}\right)\N_i} 
\\\\ \ds
~\leq~\left({\lambda_l\over \lambda_m}\right)^{\sum_{i\in\V\cap\{1,\ldots,m-1\}}\left({q-1\over q}\right)\N_i-\alpha_i-\left(1-{1\over r}\right)\left({p-1\over p}\right)\N_i+\gamma} 
\left({\lambda_l\over \lambda_m}\right)^{\sum_{i\in\V\cap\{m,\ldots,l\}}   \left({q-1\over q}\right)\N_i-\alpha_i-\left(1-{1\over r}\right)\left({p-1\over p}\right)\N_i }
\\\\ \ds~~~~~~~
\left({\lambda_l\over \lambda_m}\right)^{\sum_{i\in\V^c\cap\{1,\ldots,m-1\}}   \left({q-1\over q}\right)\N_i-\alpha_i-\left(1-{1\over r}\right)\left({p-1\over p}\right)\N_i }
\prod_{i\in\V^c\cap\{m,\ldots,l\}} \left({\lambda_l\over \lambda_i}\right)^{\left({q-1\over q}\right)\N_i-\alpha_i-\left(1-{1\over r}\right)\left({p-1\over p}\right)\N_i} 
\\\\ \ds
~\leq~\left({\lambda_l\over \lambda_m}\right)^{\sum_{i\in\V\cap\{1,\ldots,l\}}   \left({q-1\over q}\right)\N_i-\alpha_i-\left(1-{1\over r}\right)\left({p-1\over p}\right)\N_i +\gamma}.
\end{array}
\eeq
By bringing the estimates in (\ref{Case3 Indices Est1-})-(\ref{Case3 Indices Est2-}) back to (\ref{Decay Case3  Est}), we find
 \bel{Case3 Indices Result-}
 \begin{array}{lr}\ds
 \left(\lambda_m\right)^{{1\over qr}\sum_{i=1}^m \N_i-\gamma}\left(\lambda_l\right)^{\left({p-1\over pr}\right)\sum_{i=1}^l\N_i-\delta}\prod_{i=1}^n\left(\lambda_i\right)^{\alpha_i-\N_i\left({1\over p}-{1\over q}\right)}\prod_{i=1}^m\left({1\over \lambda_i}\right)^{\N_i\over qr} \prod_{i=1}^l \left({1\over \lambda_i}\right)^{\left({p-1\over pr}\right)\N_i}
\\\\ \ds
~=~\prod_{i=l+1}^n\left({\lambda_i\over \lambda_l}\right)^{\alpha_i-\N_i\left({1\over p}-{1\over q}\right)}\prod_{i=1}^{m} \left({\lambda_m\over \lambda_i}\right)^{{\N_i\over r}-\alpha_i+\left(1-{1\over r}\right)\left({1\over p}-{1\over q}\right)\N_i}
\\\\ \ds~~~~~~~
\left({\lambda_l\over \lambda_m}\right)^{\sum_{i=1}^{m-1}\left({q-1\over q}\right)\N_i-\alpha_i-\left(1-{1\over r}\right)\left({p-1\over p}\right)\N_i+\gamma}
\prod_{i=m}^{l} \left({\lambda_l\over \lambda_i}\right)^{\left({q-1\over q}\right)\N_i-\alpha_i-\left(1-{1\over r}\right)\left({p-1\over p}\right)\N_i} 
\\\\ \ds
~\leq~
\left({\lambda_l\over \lambda_m}\right)^{\sum_{i\in\V\cap\{1,\ldots,l\}}   \left({q-1\over q}\right)\N_i-\alpha_i-\left(1-{1\over r}\right)\left({p-1\over p}\right)\N_i +\gamma}
\prod_{i=1}^{m} \left({\lambda_m\over \lambda_i}\right)^{{\N_i\over r}-\alpha_i+\left(1-{1\over r}\right)\left({1\over p}-{1\over q}\right)\N_i}
\prod_{i=l+1}^n\left({\lambda_i\over \lambda_l}\right)^{\alpha_i-\N_i\left({1\over p}-{1\over q}\right)}.
\end{array}
\eeq
Recall that $\gamma>0$ satisfies the second strict inequality in (\ref{Constraint Case Three}). From (\ref{strict subbalance}) and (\ref{alpha}), we also have $\left({1\over p}-{1\over q}\right)\N_i<\alpha_i<\N_i$ for every $i=1,2,\ldots,n$. Define implicitly $0\leq\vartheta_1\leq\vartheta_2\leq1$ by letting $\lambda_m=(\lambda_n)^{\vartheta_1}$ and $\lambda_l=(\lambda_n)^{\vartheta_2}$.

For $r$ sufficiently close to $1$, we have
\bel{Case3- Indices Est}
\begin{array}{lr}\ds
~~~~~~~\vartheta_1\left[ {\N_1\over r}-\alpha_1+\left(1-{1\over r}\right)\left({1\over p}-{1\over q}\right)\N_1 \right]
~+~(1-\vartheta_2)\left( \alpha_n-\N_n\left({1\over p}-{1\over q}  \right)   \right)  
\\\\ \ds~~~~~~~
~+~
(\vartheta_2-\vartheta_1)\left[\sum_{i\in\V\cap\{1,\ldots,l\}}   \left({q-1\over q}\right)\N_i-\alpha_i-\left(1-{1\over r}\right)\left({p-1\over p}\right)\N_i +\gamma\right]
~>~0.
\end{array}
\eeq
The estimate in (\ref{Case3- Indices Est}) implies that
(\ref{Case3 Indices Result-}) is bounded by a constant multiple of $(\lambda_n)^\ve$ for some $\ve=\ve(p,q,r,\alpha,\gamma,\delta,n,\N)>0$.
\endproof

\section{One-weight inequality on product spaces}
\setcounter{equation}{0}
Let  $\omega(x)=|x|^{-\gamma},\sigma(x)=|x|^\delta, \gamma,\delta\in\R$. Consider $\gamma+\delta=0$ so that  $\omega=\sigma$. 

From (\ref{Formula}) and (\ref{subbalance}),   we must have
\bel{balance i}
{\alpha_i\over\N_i}~=~{1\over p}-{1\over q},\qquad i=1,2,\ldots,n.
\eeq
Write $x=(x_i,x_i^\dagger)\in\R^{\N_i}\times\R^{\N-\N_i}$  and $\Q^\dagger_i=\bigotimes_{i\neq j}\Q_i$ for every $i=1,2,\ldots,n$. 

Let $\Q^\dagger_i$ shrink to $x^\dagger_i$ in  (\ref{A-Characteristic}). By applying the Lebesgue Differentiation Theorem,  we have
 \bel{A-subspace i}
\begin{array}{lr}\ds
\left\{{1\over|\Q_i|}\int_{\Q_i} \omega^q(x_i,x^\dagger_i)dx_i\right\}^{1\over q}\left\{{1\over|\Q_i|}\int_{\Q_i} \left({1\over \omega}\right)^{p\over p-1}(x_i,x^\dagger_i)dx_i\right\}^{p-1\over p}~<~\infty
\end{array}
\eeq
for every  $\Q_i\subset\R^{\N_i}$ and $a\cdot e$ $x^\dagger_i\in\R^{\N-\N_i},~i=1,2,\ldots,n$.  

Observe that (\ref{balance i})-(\ref{A-subspace i}) are sufficient conditions of the Muckenhoupt-Wheeden Theorem \cite{Muckenhoupt-Wheeden} which implies  
\bel{one weight i}
\begin{array}{lr}\ds
\left\{\int_{\R^{\N_i}} \left\{\int_{\R^{\N_i}} f(y_i,x_i^\dagger) \left({1\over |x_i-y_i|}\right)^{\N_i-\alpha_i} dy_i\right\}^q\omega^q(x_i,x_i^\dagger)dx_i\right\}^{1\over q}
\\\\ \ds~~~~~~~~~~~~~~~~~~~~~~~~~~~~~~~~~
~\leq~\C_{p~q~\N_i~\omega} \left\{\int_{\R^{\N_i}} \Big(f\omega\Big)^p(x_i,x_i^\dagger)dx_i \right\}^{1\over p}
\end{array}
\eeq
 for $1<p<q<\infty$ and $a.e$ $x_i^\dagger\in\R^{\N-\N_i},~i=1,2,\ldots,n$. 
 
 By using (\ref{one weight i}), we   have 
 \bel{One-Weight}
\begin{array}{lr}\ds
\left\{\int_{\R^\N}   \Big(\omega\I_\alpha f\Big)^q(x)dx \right\}^{1\over q}
~=~\left\{\int_{\R^\N} \left\{\int_{\R^\N} f(y)\prod_{i=1}^n \left({1\over |x_i-y_i|}\right)^{\N_i-\alpha_i}dy\right\}^q \omega^q(x)dx\right\}^{1\over q}
\\\\ \ds
~\leq~\C_{p~q~\N_i~\omega}\left\{ \int_{\R^{\N-\N_i}}\left\{ \int_{\R^{\N_i}}\left\{\int_{\R^{\N-\N_i}} f(x_i,y^\dagger_i)\prod_{j\neq i} \left({1\over |x_j-y_j|}\right)^{\N_j-\alpha_j} dy^\dagger_i\right\}^p \omega^p(x_i,x^\dagger_i)dx_i\right\}^{q\over p} dx^\dagger_i\right\}^{1\over q}
\\\\ \ds
~\leq~\C_{p~q~\N_i~\omega}\left\{ \int_{\R^{\N_i}}\left\{ \int_{\R^{\N-\N_i}}\left\{\int_{\R^{\N-\N_i}} f(x_i,y^\dagger_i)\prod_{j\neq i} \left({1\over |x_j-y_j|}\right)^{\N_j-\alpha_j} dy^\dagger_i\right\}^q \omega^q(x_i,x^\dagger_i)dx^\dagger_i\right\}^{p\over q} dx_i\right\}^{1\over p}
\\ \ds~~~~~~~~~~~~~~~~~~~~~~~~~~~~~~~~~~~~~~~~~\vdots~~~~~~~~~~~~~~~~~~~~~~~~~~~~~~~~~~~~~~~~~~~~~~~
\hbox{\small{by Minkowski integral inequality}}
\\ \ds 
~\leq~\C_{p~q~n~\N~\omega}\left\{\int_{\R^\N} \Big(f\omega\Big)^p (x)dx\right\}^{1\over p},\qquad 1<p<q<\infty.
\end{array}
\eeq

\section{Proof of Theorem A*}
\setcounter{equation}{0}
Let $\{1,2,\ldots n\}=\setI\cup\setJ$ where
\bel{IJ}
\setI~=~\left\{i\in\{1,2,\ldots,n\}~\colon~{\alpha_i\over\N_i}
~=~{1\over p}-{1\over q}\right\},\qquad
\setJ~=~\left\{i\in\{1,2,\ldots,n\}~\colon~{\alpha_i\over\N_i}
~>~{1\over p}-{1\over q}\right\}.
\eeq
Define
\bel{setting IJ}
\begin{array}{lr}\ds
\alpha_\setI~=~\sum_{i\in\setI}\alpha_i,~~~~\Q_\setI=\bigotimes_{i\in\setI}\Q_i,~~~~\R^{\N_\setI}~=~\bigotimes_{i\in\setI}\R^{\N_i},
\\\\ \ds
\alpha_\setJ~=~\sum_{i\in\setJ}\alpha_i,~~~~\Q_\setJ~=~\bigotimes_{i\in\setJ}\Q_i,~~~~ \R^{\N_\setJ}~=~\bigotimes_{i\in\setJ}\R^{\N_i}.
 \end{array}
\eeq  
We write $x=(x_\setI,x_\setJ)\in\R^{\N_\setI}\times\R^{\N_\setJ}$ and denote the cardinality of $\setI$ and $\setJ$ by  $|\setI|$ and $|\setJ|$. 

Suppose  $\omega(x)=|x|^{-\gamma},\sigma(x)=|x|^\delta,~\gamma,\delta\in\R$ satisfy the Muckenhoupt characteristic  (\ref{A-Characteristic}). 
Consider $\Q_i$ centered on the origin of $\R^{\N_i}$ for every $i\in\setI$.   
Let $\Q_i, i\in\setI$ shrink to the origin.
By applying Lebesgue Differentiation Theorem, we have
\bel{A-subspace}
\begin{array}{lr}\ds
\sup_{\Q_\setJ\subset\R^{\N_\setJ}}~\prod_{i\in\setJ} |\Q_i|^{{\alpha_i\over \N_i}-\left({1\over p}-{1\over q}\right)}\left\{{1\over|\Q_\setJ|}\int_{\Q_\setJ} \left({1\over|x_\setJ|}\right)^{q\gamma}dx_\setJ\right\}^{1\over q}\left\{{1\over|\Q_\setJ|}\int_{\Q_\setJ} \left({1\over|x_\setJ|}\right)^{\delta p\over p-1}dx_\setJ\right\}^{p-1\over p}
\\\\ \ds
~\leq~\sup_{\Q\subset\R^\N}~\prod_{i\in\setJ} |\Q_i|^{{\alpha_i\over \N_i}-\left({1\over p}-{1\over q}\right)}\left\{{1\over|\Q|}\int_{\Q} \left({1\over|x_\setJ|}\right)^{q\gamma}dx_\setJ\right\}^{1\over q}\left\{{1\over|\Q|}\int_{\Q} \left({1\over|x_\setJ|}\right)^{\delta p\over p-1}dx_\setJ\right\}^{p-1\over p}~<~\infty.
\end{array}
\eeq
The boundedness of $\A_{pq}^\alpha\left(|x|^{-\gamma},|x|^\delta\right)$ requires $\gamma q<\N_\setJ$ and $\delta\left({p\over p-1}\right)<\N_\setJ$.

\begin{prop} Let $\omega(x_\setJ)=|x_\setJ|^{-\gamma},\sigma(x_\setJ)=|x_\setJ|^\delta$ for $ \gamma,\delta\in\R$ satisfying (\ref{A-subspace}). For $a.e ~x_\setI\in\R^{\N_\setI}$, we have
\bel{Norm Ineq Case3}
\begin{array}{lr}\ds
\left\{\int_{\R^{\N_\setJ}} \left\{ \int_{\R^{\N_\setJ}} f(x_\setI,y_\setJ)\prod_{i\in\setJ}\left({1\over |x_i-y_i|}\right)^{\N_i-\alpha_i}dy_\setJ\right\}^q\omega^q(x_\setJ) dx_\setJ\right\}^{1\over q}
\\\\ \ds
~\leq~\C_{p~q~\alpha_\setJ~\gamma~\delta~|\setJ|~\N_\setJ}\left\{\int_{\R^{\N_\setJ}} \Big(f(x_\setI,x_\setJ)\Big)^p\sigma^p(x_\setJ)dx_\setJ\right\}^{1\over p},\qquad 1~<~p~\leq ~q~<~\infty.
\end{array}
\eeq
\end{prop}

{\bf Proof:} From section 4, we have the Muckenhoupt characteristic (\ref{A-subspace}) implying $\gamma, \delta$ to satisfy (\ref{local integrability})-(\ref{Constraint Case Three}) with $\alpha, n, \N$ replaced by $\alpha_\setJ$, $|\setJ|$, $\N_\setJ$ respectively. 
Suppose $|\setJ|=1$. {\bf Theorem A} by Stein and Weiss \cite{Stein-Weiss} shows that these constraints are   sufficient conditions to imply  (\ref{Norm Ineq Case3}). 

Consider $|\setJ|\ge2$. By applying {\bf Principal Lemma} in the beginning of section 5, we have $\omega(x_\setJ)=|x_\setJ|^{-\gamma},\sigma(x_\setJ)=|x_\setJ|^\delta$ satisfying the decay estimate  (\ref{ratio Q})-(\ref{Eccentricity Decay}) for every  $\Q_\setJ\subset\R^{\N_\setJ}$  where $\alpha_i>\N_i\left({1\over p}-{1\over q}\right),~i\in\setJ$.

Let $\t$ denote the $|\setJ|$-tuple $\left(2^{-t_1},2^{-t_2},\ldots,2^{-t_{|\setJ|}}\right)$.
We have
\bel{Summability Case3}
\sum_\t \A_{pqs}^{\alpha_\setJ} \left(\t~\colon |x_\setJ|^{-\gamma},|x_\setJ|^\delta\right)~<~\infty
\eeq
as required in (\ref{Summability}) for every $0<s<1$.
\endproof
\v
Let $\gamma>0,\delta>0$ satisfy  (\ref{local integrability})-(\ref{Formula}) and (\ref{Constraint Case Three}). In particular,  we have
\bel{Compara>>}
 \omega(x)=|x|^{-\gamma}~\leq~|x_\setJ|^{-\gamma}=\omega(x_\setJ),\qquad  \sigma(x_\setJ)=|x_\setJ|^\delta~\leq~|x|^\delta=\sigma(x).
 \eeq
From (\ref{Compara>>}), we have
\bel{Reduction Case Three}
\begin{array}{lr}\ds
\left\{\int_{\R^\N} \Big(\omega \I_\alpha f\Big)^q(x)dx\right\}^{1\over q}
~\leq~\left\{\int_{\R^\N} \left\{ \int_{\R^\N}f(y) \prod_{i=1}^n\left({1\over |x_i-y_i|}\right)^{\N_i-\alpha_i} dy\right\}^q \omega^q\left(x_\setJ\right)dx\right\}^{1\over q}
\\\\ \ds
~\leq~\C_{p~q~\alpha_\setJ~\gamma~\delta~|\setJ|~\N_\setJ}\left\{\int_{\R^{\N_\setI}} \left\{\int_{\R^{\N_\setJ}} \left\{\int_{\R^{\N_\setI}}f\left(y_\setI,x_\setJ\right)\prod_{i\in\setI}\left({1\over |x_i-y_i|}\right)^{\N_i-\alpha_i}  dy_\setI\right\}^p \sigma^p\left(x_\setJ\right)dx_\setJ\right\}^{q\over p} dx_\setI  \right\}^{1\over q}
\\ \ds~~~~~~~~~~~~~~~~~~~~~~~~~~~~~~~~~~~~~~~~~~~~~~~~~~~~~~~~~~~~~~~~~~~~~~~~~~~~~~~~~~~~~~~~~~~~~~~~~~~~~~~~~~~~~~~~~~~~~~~~~~~~~~~
\hbox{\small{by {\bf Proposition 7.1}}}
\\ \ds
~\leq~\C_{p~q~\alpha_\setJ~\gamma~\delta~|\setJ|~\N_\setJ}\left\{\int_{\R^{\N_\setJ}} \left\{\int_{\R^{\N_\setI}} \left\{\int_{\R^{\N_\setI}}f\left(y_\setI,x_\setJ\right)\prod_{i\in\setI}\left({1\over |x_i-y_i|}\right)^{\N_i-\alpha_i}  dy_\setI\right\}^q dx_\setI\right\}^{p\over q}\sigma^p\left(x_\setJ\right) dx_\setJ  \right\}^{1\over p}
\\ \ds~~~~~~~~~~~~~~~~~~~~~~~~~~~~~~~~~~~~~~~~~~~~~~~~~~~~~~~~~~~~~~~~~~~~~~~~~~~~~~~~~~~~~~~~~~~~~~~~~~~~~~~
\hbox{\small{by Minkowski integral inequality}}
\\ \ds
~\leq~\C_{p~q~\alpha~\gamma~\delta~n~\N}\left\{\iint_{\R^{\N_\setI}\times\R^{\N_\setJ}} \Big(f\left(x_\setI,x_\setJ\right)\Big)^p\sigma^p\left(x_\setJ\right) dx_\setI dx_\setJ\right\}^{1\over p}\qquad\hbox{\small{by (\ref{A-subspace i})-(\ref{One-Weight}) }}
\\\\ \ds
~\leq~\C_{p~q~\alpha~\gamma~\delta~n~\N}\left\{\int_{\R^\N} \Big(f\sigma\Big)^p(x)dx\right\}^{1\over p},\qquad 1<p\leq q<\infty.  
\end{array}
\eeq
Consider $\gamma\ge0,\delta\leq0$ satisfying (\ref{local integrability})-(\ref{Constraint Case One}) or $\gamma\leq0,\delta\ge0$ satisfying (\ref{local integrability})-(\ref{Formula}) and (\ref{Constraint Case Two}). 
Note that it is suffice to study one of these two cases because  
 $\I_\alpha$ is  self-adjoint  and
\bel{duality}
\left\|\omega \I_\alpha f\right\|_{\L^q\left(\R^\N\right)}~\lesssim~\left\| f\sigma\right\|_{\L^p\left(\R^\N\right)}\qquad\hbox{if and only if}\qquad \left\|\sigma^{-1} \I_\alpha g\right\|_{\L^{p\over p-1}\left(\R^\N\right)}~\lesssim~\left\| g\omega^{-1}\right\|_{\L^{q\over q-1}\left(\R^\N\right)}.
\eeq
Let $\gamma\ge0,\delta\leq0$. Suppose that $f$ is supported in the region where $|x_\setI|\leq|x_\setJ|$. 
By using (\ref{Compara>>}) and carrying out the same estimate (\ref{Reduction Case Three}), we have
\bel{Reduction Case One easy}
\begin{array}{lr}\ds
\left\{\int_{\R^\N} \left\{\int_{ |y_\setI|\leq|y_\setJ| } f(y)\prod_{i=1}^n\left({1\over |x_i-y_i|}\right)^{\N_i-\alpha_i} dy\right\}^q  \omega^q(x)dx\right\}^{1\over q}
\\\\ \ds
~\leq~\left\{\iint_{\R^{\N_\setI}\times\R^{\N_\setJ}} \left\{\iint_{|y_\setI|\leq|y_\setJ|} f\left(y_\setI,y_\setJ\right)\prod_{i\in\setI\cup\setJ}\left({1\over |x_i-y_i|}\right)^{\N_i-\alpha_i} dy_\setI dy_\setJ\right\}^q  \omega^q(x_\setJ)dx\right\}^{1\over q}
\\\\ \ds
~\leq~\C_{p~q~\alpha~\gamma~\delta~n~\N}\left\{\iint_{|x_\setI|\leq|x_\setJ|} \Big( f(x_\setI,x_\setJ)\Big)^p \sigma^p(x_\setJ) dx_\setI dx_\setJ\right\}^{1\over p}
\\\\ \ds
~\leq~\C_{p~q~\alpha~\gamma~\delta~n~\N}\left\{\int_{\R^\N} \Big(f\sigma\Big)^p(x)dx\right\}^{1\over p}.\end{array}
\eeq
 The last inequality holds in (\ref{Reduction Case One easy}) because $\sigma(x_\setJ)=|x_\setJ|^\delta\approx|x|^\delta=\sigma(x)$ for $|x_\setI|\leq|x_\setJ|$.

On the other hand, suppose $f$ supported in the region $|x_\setI|>|x_\setJ|$.
Recall that $\gamma\ge0,\delta\leq0$ satisfy (\ref{local integrability})-(\ref{Constraint Case One}). 
In particular,  
 \bel{Nece Case One}
\gamma+\delta~=~\sum_{i=1}^n \alpha_i-\N_i\left({1\over p}-{1\over q}\right)~=~\sum_{i\in\setJ}\alpha_i-\N_i\left({1\over p}-{1\over q}\right)\qquad\hbox{\small{by (\ref{IJ})}}
\eeq
and
\bel{Nece Case One'}
\alpha_i-{\N_i\over p}~<~\delta~\leq~0~~~~ \hbox{for every}~~~~ i\in\{1,2,\ldots,n\}=\setI\cup\setJ.
\eeq
By putting together (\ref{Nece Case One}) and (\ref{Nece Case One'}), we find
\bel{constraint case one rho}
0~\leq~\gamma+\delta~=~\sum_{i\in\setJ}\alpha_i-\N_i\left({1\over p}-{1\over q}\right)~<~{\N_\setJ\over q},\qquad 0~<~\N_\setJ\left({p-1\over p}\right).
\eeq
\begin{prop} Let $\rho\left(x_\setJ\right)=|x_\setJ|^{-(\gamma+\delta)},~\eta\left(x_\setJ\right)\equiv1$.  For $a.e$ $x_\setI\in\R^{\N_\setI}$, we have
\bel{Norm Ineq Case1}
\begin{array}{lr}\ds
\left\{\int_{\R^{\N_\setJ}} \left\{ \int_{\R^{\N_\setJ}} f(x_\setI,y_\setJ)\prod_{i\in\setJ}\left({1\over |x_i-y_i|}\right)^{\N_i-\alpha_i}dy_\setJ\right\}^q\rho^q(x_\setJ) dx_\setJ\right\}^{1\over q}
\\\\ \ds
~\leq~\C_{p~q~\alpha_\setJ~\gamma~\delta~|\setJ|~\N_\setJ}\left\{\int_{\R^{\N_\setJ}} \Big(f(x_\setI,x_\setJ)\Big)^p\eta^p(x_\setJ)dx_\setJ\right\}^{1\over p},\qquad 1~<~p~\leq ~q~<~\infty.
\end{array}
\eeq
\end{prop}  
{\bf Proof:} Observe that (\ref{Nece Case One})-(\ref{constraint case one rho}) imply the constraints in (\ref{local integrability})-(\ref{Constraint Case One})  with $\gamma,\delta,\alpha, n, \N$ replaced by $\gamma+\delta, 0, \alpha_\setJ, |\setJ|,\N_\setJ$ respectively.   

Suppose $|\setJ|=1$, From  {\bf Theorem A}, it follows that (\ref{Nece Case One})-(\ref{constraint case one rho}) are sufficient conditions to imply (\ref{Norm Ineq Case1}). 

Consider $|\setJ|\ge2$.  By applying {\bf Principal Lemma},  $\rho(x_\setJ)=|x_\setJ|^{-(\gamma+\delta)}, \eta(x_\setJ)\equiv1$ satisfy the decay estimate in (\ref{ratio Q})-(\ref{Eccentricity Decay}) for  every $\Q_\setJ\subset\R^{\N_\setJ}$ where $\alpha_i>\N_i\left({1\over p}-{1\over q}\right),~i\in\setJ$. 

Let $\t$ denote the $|\setJ|$-tuple $\left(2^{-t_1},2^{-t_2},\ldots,2^{-t_{|\setJ|}}\right)$. We have
\bel{Summability Case1}
\sum_\t \A_{pqs}^{\alpha_\setJ} \left(\t~\colon |x_\setJ|^{-(\gamma+\delta)},1\right)~<~\infty
\eeq
as required in (\ref{Summability}) for every $0<s<1$. \endproof

\begin{prop} Let $\omega(x_\setI)=\sigma\left(x_\setI\right)=|x_\setI|^{\delta}$. For $a.e$ $x_\setJ\in\R^{\N_\setJ}$, 
we have
\bel{Norm Ineq One-Weight}
\begin{array}{lr}\ds
\left\{\int_{\R^{\N_\setI}} \left\{ \int_{\R^{\N_\setI}} f(y_\setI,x_\setJ)\prod_{i\in\setI}\left({1\over |x_i-y_i|}\right)^{\N_i-\alpha_i}dy_\setI\right\}^q\omega^q(x_\setI) dx_\setI\right\}^{1\over q}
\\\\ \ds
~\leq~\C_{p~q~\alpha_\setJ~\gamma~\delta~|\setI|~\N_\setI}\left\{\int_{\R^{\N_\setI}} \Big(f(x_\setI,x_\setJ)\Big)^p\omega^p(x_\setI)dx_\setI\right\}^{1\over p},\qquad 1~<~p~< ~q~<~\infty.
\end{array}
\eeq
\end{prop}
{\bf Proof:} Recall (\ref{IJ}) and (\ref{Nece Case One'}). We have
\bel{constraint case one I}
\begin{array}{lr}\ds
-\delta+\delta~=~0~=~\sum_{i\in\setI}\alpha_i-\N_i\left({1\over p}-{1\over q}\right),
\qquad
-\delta~<~{\N_i\over p}-\alpha~=~{\N_i\over q}~~~\hbox{for}~~~ i\in\setI.
\end{array}
\eeq
Note that the constraints in (\ref{constraint case one I}) are sufficient conditions for {\bf Theorem A} on every subspace $\R^{\N_i}, i\in\setI$. The norm inequality  (\ref{Norm Ineq One-Weight}) can be obtained by  following the iteration argument given in section 6.
\endproof

Let $\rho(x_\setJ)=|x_\setJ|^{-(\gamma+\delta)}$ and  $\sigma(x_\setI)=|x_\setI|^\delta$ where $\gamma+\delta\ge0$ and $\delta\leq0$. It is clear that
\bel{Compara><}
\omega(x)~=~|x|^{-\gamma}~\leq~\rho\left(x_\setJ\right)\sigma\left(x_\setI\right).
\eeq
We have
\bel{Reduction Case One}
\begin{array}{lr}\ds
\left\{\int_{\R^\N} \left\{ \int_{|y_\setI|>|y_\setJ|}f(y)\prod_{i=1}^n\left({1\over |x_i-y_i|}\right)^{\N_i-\alpha_i} dy\right\}^q \omega^q(x)dx\right\}^{1\over q}
\\\\ \ds
~\leq~\left\{\int_{\R^\N} \left\{ \int_{|y_\setI|>|y_\setJ|}f(y) \prod_{i=1}^n\left({1\over |x_i-y_i|}\right)^{\N_i-\alpha_i}dy\right\}^q \rho^q\left(x_\setJ\right)\sigma^q\left(x_\setI\right)dx\right\}^{1\over q}\qquad \hbox{\small{by (\ref{Compara><})}}
\\\\ \ds
~=~\left\{\iint_{\R^{\N_\setI}\times\R^{\N_\setJ}} \left\{\iint_{|y_\setI|>|y_\setJ|} f(y_\setI,y_\setJ)\prod_{i\in\setI\cup\setJ}\left({1\over |x_i-y_i|}\right)^{\N_i-\alpha_i}dy_\setI dy_\setJ\right\}^q  \rho^q\left(x_\setJ\right)\sigma^q\left(x_\setI\right)dx_\setI dx_\setJ\right\}^{1\over q}

\\\\ \ds
~\leq~\C_{p~q~\alpha_\setJ~\gamma~\delta~|\setJ|~\N_\setJ}\left\{\int_{\R^{\N_\setI}} \left\{\int_{\R^{\N_\setJ}} \left\{\int_{|y_\setI|>|x_\setJ|} f (y_\setI,x_\setJ)\prod_{i\in\setI}\left({1\over |x_i-y_i|}\right)^{\N_i-\alpha_i}  dy_\setI\right\}^p dx_\setJ\right\}^{q\over p} \sigma^q\left(x_\setI\right)dx_\setI  \right\}^{1\over q}
\\ \ds~~~~~~~~~~~~~~~~~~~~~~~~~~~~~~~~~~~~~~~~~~~~~~~~~~~~~~~~~~~~~~~~~~~~~~~~~~~~~~~~~~~~~~~~~~~~~~~~~~~~~~~~~~~~~~~~~~~~~~~~~~~~\hbox{\small{by {\bf Proposition 7.2}}}
\\ \ds
~\leq~\C_{p~q~\alpha_\setJ~\gamma~\delta~|\setJ|~\N_\setJ}\left\{\int_{\R^{\N_\setJ}} \left\{\int_{\R^{\N_\setI}} \left\{\int_{|y_\setI|>|x_\setJ|}f(y_\setI,x_\setJ) \prod_{i\in\setI}\left({1\over |x_i-y_i|}\right)^{\N_i-\alpha_i}  dy_\setI\right\}^q \sigma^q\left(x_\setI\right)dx_\setI\right\}^{p\over q} dx_\setJ  \right\}^{1\over p}
\\ \ds~~~~~~~~~~~~~~~~~~~~~~~~~~~~~~~~~~~~~~~~~~~~~~~~~~~~~~~~~~~~~~~~~~~~~~~~~~~~~~~~~~~~~~~~~~~~~~~~~~~~~~~~~
\hbox{by Minkowski integral inequality}
\\ \ds
~\leq~\C_{p~q~\alpha~\gamma~\delta~n~\N}\left\{\iint_{|x_\setI|>|x_\setJ|}\Big(f(x_\setI,x_\setJ)\Big)^p   \sigma^p\left(x_\setI\right) dx_\setI dx_\setJ\right\}^{1\over p}
\qquad \hbox{\small{by {\bf Proposition 7.3}}}
\\\\ \ds
~\leq~\C_{p~q~\alpha~\gamma~\delta~n~\N}\left\{\int_{\R^\N} \Big(f\sigma\Big)^p(x)dx\right\}^{1\over p},\qquad1<p\leq q<\infty.\end{array}
\eeq
The last inequality holds  because $\sigma(x_\setI)\approx\sigma(x)$ for $|x_\setI|>|x_\setJ|$.
\v

{\bf Acknowledgement:

~~~~~~~I am deeply grateful to my advisor Elias M. Stein for those stimulating talks and unforgettable lectures.}

\end{document}